\documentclass[reqno]{amsart}

\usepackage[margin=1in]{geometry}

\usepackage{graphicx,epstopdf,epsfig}
\usepackage{amsfonts,epsfig,fancyhdr,graphics, hyperref,amsmath,amssymb}
\usepackage[inline]{enumitem}
\usepackage{mathtools}
\usepackage{tikz,pgfplots}
\usepgfplotslibrary{fillbetween}

\usepackage{amsthm}
\usepackage{cleveref}

\makeatletter
\g@addto@macro{\endabstract}{\@setabstract}
\newcommand{\authorfootnotes}{\renewcommand\thefootnote{\@fnsymbol\c@footnote}}%
\makeatother

\theoremstyle{plain}
\newtheorem{theorem}{Theorem}[section]
\newtheorem{corollary}[theorem]{Corollary}
\newtheorem{proposition}[theorem]{Proposition}
\newtheorem{lemma}[theorem]{Lemma}

\theoremstyle{definition}
\newtheorem{definition}[theorem]{Definition}
\newtheorem{remark}[theorem]{Remark}
\newtheorem{example}[theorem]{Example}

\newcommand{\bd}{\mathbf d}
\newcommand{\be}{\mathbf e}

\newcommand{\cP}{\mathcal P}
\newcommand{\cL}{\mathcal L}

\title{Kemeny's constant and enumerating Braess edges in trees}

\begin{document}

\maketitle

\authorfootnotes
Jihyeug Jang\textsuperscript{1}, Mark Kempton\textsuperscript{2},
Sooyeong Kim\footnote{Contact: kimswim@yorku.ca}\textsuperscript{3}, Adam Knudson\textsuperscript{2}, Neal Madras\textsuperscript{3} and
Minho Song\textsuperscript{1} \par \bigskip

\textsuperscript{1}Department of Mathematics, Sungkyunkwan University, Suwon, 16419, South Korea \par
\textsuperscript{2}Department of Mathematics, Brigham Young University, Provo UT, USA\par 
\textsuperscript{3}Department of Mathematics and Statistics, York University, 4700 Keele Street, Toronto, Canada\par \bigskip

\begin{abstract}
    We study the problem of enumerating Braess edges for Kemeny's constant in trees.  We obtain bounds and asympotic results for the number of Braess edges in some families of trees.
\end{abstract}



\section{Introduction}
An important parameter arising from the study of Markov chains that has received increased attention in graph theory in recent years is Kemeny's constant.  Given a discrete time finite state irreducible Markov chain with stationary distribution $\pi$, we define \emph{Kemeny's constant} of the Markov chain by
\[
\kappa = \sum_{j=1}^nm_{ij}\pi(j)
\]
where $m_{ij}$ denotes the \emph{mean first passage time} from $i$ to $j$, which is the expected number of steps for the Markov chain to go from state $i$ to state $j$.  Surprisingly, this quantity is independent of the state $i$ \cite{KemenySnell}. In this paper, we will be studying Kemeny's constant of a random walk on a graph. As a graph parameter, Kemeny's constant provides a measure of how hard it is for a random walker to move around the graph, and gives a way of measuring how well connected the graph is---a small Kemeny's constant corresponds to a well-connected graph in which it is easy to move around, while a large Kemeny's constant corresponds to a poorly connected graph, possibly with bottlenecks. For a graph $G$, we will denote Kemeny's constant for the random walk on $G$ as $\kappa(G)$, and we will simply call it Kemeny's constant of $G$. Kemeny's constant has received considerable attention recently in the study of graphs and networks, and has many applications to a variety of subjects, including the study of disease spread \cite{kahn2014mixing,kim2023effectiveness,ruhi2015sirs,van2008virus,yilmaz2020kemeny}, network algorithms \cite{dudkina2021node,fouss2016algorithms,levene2002kemeny}, robotic surveillance \cite{patel2015robotic}, and many others. 

One interesting phenomenon studied in connection with Kemeny's constant is \emph{Braess' Paradox}.  Loosely speaking, Braess' paradox occurs when the deletion of an edge from (or insertion of a new edge to) a graph affects the overall connectivity in a counterintuitive way.   Originally, Braess \cite{braess2005paradox} studied road networks, and observed instances in road network models in which the closure of a road would actually improve overall traffic flow in a city.  Intuitively, one expects that adding an edge to a graph will improve the graph's overall connectivity.  When this is not the case, we say that Braess' paradox occurs.  Braess' paradox in traffic networks has been widely studied since this initial observation; see for instance \cite{ding2012traffic,hagstrom2001characterizing,rapoport2009choice} and references therein.

The study of Braess' paradox in the context of Kemeny's constant was first introduced by Kirkland and Zeng in \cite{kirkland2016kemeny}.  We say a pair of nonadjacent vertices $\{ u,v \}$ in a graph $G$ is a \emph{Braess edge} of $G$ if adding the edge $\{ u,v \}$ to $G$ increases Kemeny's constant.  Kirkland and Zeng proved that if $G$ is a tree with two pendent twin vertices (degree one vertices adjacent to the same vertex), then the endpoints of those pendent vertices constitute a Braess edge.  This fact was generalized by Ciardo \cite{ciardo2022kemeny} to arbitrary connected graphs with a pair of pendent twins.  Intuitively, when the edge between the twin pendent vertices is added, then there is some probability that a random walker will get stuck walking between those two vertices for some time, thus increasing travel times to other parts of the graph.  Work in \cite{faught20221} looked more generally at graphs with a cut vertex (a vertex whose deletion disconnects the graph) and found that Braess edges occur frequently in such graphs. In particular, graphs with a cut vertex whose deletion gives rise to two paths were studied by Kim \cite{kim2022families}, and it was discovered that the eccentricity of the cut vertex is related to the formation of the Braess edge joining the end-vertices of each path. Furthermore, work by Hu and Kirkland \cite{hu2019complete} established equivalent conditions for complete multipartite graphs and complete split graphs to have every non-edge as a Braess edge.  Recently, work from \cite{kirkland2023edge} studied the overall effect that adding an edge to a graph can have on Kemeny's constant.

The purpose of this paper is to study the occurrence of Braess edges in trees in greater depth.  Our motivation is to understand what kind of Braess edges can occur in trees and to enumerate how many Braess edges can occur in a tree.  We are able to give partial answers to these questions for various families of trees.  

After establishing the terminology, notation, and tools we will use in Section \ref{sec:prelim}, we will first focus on paths in Section \ref{sec:path}.  Paths are of particular interest because a path on $n$ vertices is the tree on $n$ vertices with the largest Kemeny's constant, as shown in \cite{faught20221}.  We show that the number of Braess edges for a path of length $k$ is \[\frac13k\ln k - ck +o(k)\] for a constant $c\approx .548$ (see Corollary \ref{cor-asymp} below).  We further show that the addition of any edge that crosses the middle of the path cannot be a Braess edge, showing that the Braess edges in a path are concentrated towards the endpoints (see Theorems \ref{thm:nonedge crossing the center} and \ref{thm:equiv Braess on path} as well as Remark \ref{rmk:far from the center} below).  This fits our intuition: an edge connecting points that were originally far apart will tend to decrease travel times of the random walker, whereas edges close to the endpoints of a path create bottlenecks where the random walker could get stuck for a time.  

In Section \ref{sec:spider}, we study a family of trees that we call \emph{spider graphs}: the spider graph $\mathcal{S}_{a,b}$ is obtained from taking $b$ paths of length $a$ and connecting one endpoint of each at a central vertex of degree $b$.  This is motivated by thinking of stars, in which results from \cite{ciardo2022kemeny,kirkland2016kemeny} show that any pair of nonadjacent vertices is Braess. We show that the number of Braess edges for $\mathcal{S}_{a,b}$ is
$$\Theta(ab\ln(a))$$
for sufficiently large $a$.
(See Theorem~\ref{thm:spider} for lower and upper bounds.) We further show that adding an edge that crosses this central vertex cannot be a Braess edge (see Theorem~\ref{prop:nonBraess1}).  For vertices along the branches of a spider graph, we prove that in $\mathcal{S}_{a,b}$, if there is a Braess edge in a path of length $2a$ (which is $\mathcal{S}_{a,2}$), then the corresponding edge is also Braess in $\mathcal{S}_{a,b}$ (see Lemma \ref{lem:Braess on spider}). For sufficiently large $a$ and $b$, $\mathcal{S}_{a,b}$ can contain more Braess edges as well (see Theorems \ref{thm:spider} and \ref{thm:spider2}).  As a consequence of results from this section, we are able to show that there are no Braess edges in $\mathcal{S}_{2,b}$ for $b\geq2$  (see Corollary~\ref{cor:no Braess}), thus exhibiting a family of trees in which Braess' paradox does not occur for any pair of vertices.

Finally, in Section \ref{sec:broom}, we combine ideas from the previous two sections and study a family of graph that we call \emph{brooms}.  The broom $\mathcal{B}_{k,p}$ is obtained from attaching $p$ pendent vertices to the endpoint of a path of length $k$.  We determine the asympototics for the number of Braess edges in broom graphs.  In particular, there are three kinds of non-edges in a broom graph: pairs whose endpoints are both among the pendent vertices of the start at the end, those where both vertices are along the path, and those where one endpoint is on the star and the other on the path. We show that in a broom $\mathcal{B}_{k,p}$, as $k\rightarrow\infty$, there are $\binom{p}{2}$ Braess edges among the pendent vertices, $\Theta(k\ln{k})$ Braess edges along the path, and $O(k)$ Braess edges attaching a pendent vertex to the path (see Theorem \ref{thm:broom}).

\section{Preliminaries}\label{sec:prelim}

We begin with necessary terminology and notation in graph theory. Let $G$ be a graph of order $n$ with vertex set $V(G)$ and edge set $E(G)$ where $n=|V(G)|$. An edge joining vertices $v$ and $w$ of $G$ is denoted by $\{ v,w \}$. Let $m_G$ be defined as $|E(G)|$. The subgraph of $G$ \textit{induced} by a subset $S$ of $V(G)$ is the graph with vertex set $S$, where two vertices in $S$ are adjacent if and only if they are adjacent in $G$. For $v\in V(G)$, we denote by $\mathrm{deg}_G(v)$ the degree of $v$. A vertex $v$ of a graph $G$ is said to be \textit{pendent} if $\mathrm{deg}_G(v)=1$. Given a labelling of $V(G)$, we define $\mathbf{d}_G$ to be the column vector whose $i^\text{th}$ component is $\mathrm{deg}_G(v_i)$ for $1\leq i\leq n$, where $v_i$ is the $i^\text{th}$ vertex in $V(G)$. For $v,w\in V(G)$, the distance between $v$ and $w$ in $G$ is denoted by $\mathrm{dist}_G(v,w)$. We define $F_{G}$ to be the matrix given by $F_{G}=[f_{i,j}^{G}]$ where $f_{i,j}^{G}$ is the number of $2$-tree spanning forests of $G$ such that one of the two trees contains a vertex $i$ of $G$, and the other has a vertex $j$ of $G$. Note that $f_{i,i}^{G}=0$, that is, the diagonal entries of $F_G$ are zero. We denote by $\be_v$ the column vector whose component in $v^\text{th}$ position is $1$ and zeros elsewhere, and denote by $\mathbf{f}_G^v$ the $v^{\text{th}}$ column of $F_{G}$. That is, the $i^{\text{th}}$ entry of $\mathbf{f}_G^v$ is $f_{i,v}^{G}$. We define $B(G)$ to be the number of Braess edges for $G$, and denote by $\tau_G$ the number of spanning trees of $G$.

Throughout the paper, we will be using standard notation regarding asymptotic analysis of a function of the number $n$ of vertices of a graph.  In particular, we say that $f(n)$ is $O(g(n))$ (or $f(n)=O(g(n))$) if $|f(n)|\leq C|g(n)|$ for some constant $C$ not depending on $n$ for sufficiently large $n$.  We say that $f(n)$ is $o(g(n))$ if $\lim_{n\rightarrow\infty} (f(n)/g(n)) = 0$.  And finally, we say that $f(n)\sim g(n)$ if $\lim_{n\rightarrow\infty} (f(n)/g(n)) = 1$.

\begin{proposition}\cite{kim2022families}\label{Prop:dFd with a cut-vertex}
	Let $H_1$ and $H_2$ be connected graphs, and let $v_1\in V(H_1)$ and $v_2\in V(H_2)$. Assume that $G$ is obtained from $H_1$ and $H_2$ by identifying $v_1$ and $v_2$ as a vertex $v$. Suppose that $\widetilde{H}_1=H_1-v_1$ and $\widetilde{H}_2=H_2-v_2$. Then, labelling the vertices of $G$ in order of $V(\widetilde{H}_1)$, $v$, and $V(\widetilde{H}_2)$, we have:
	\begin{align*}
		\mathbf{d}_{G}^T&=[\mathbf{d}_{H_1}^T\;\mathbf{0}^T_{|V(\widetilde{H}_2)|}]+[\mathbf{0}^T_{|V(\widetilde{H}_1)|}\;\mathbf{d}_{H_2}^T],\\
		m_{G}&=m_{H_1}+m_{H_2},\\
		\tau_{G}&=\tau_{H_1}\tau_{H_2},\\
		F_{G}&=\left[\begin{array}{c|c|c}
			\tau_{H_2}F_{\widetilde{H}_1} & \tau_{H_2}\mathbf{f}_1 & \tau_{H_2}\mathbf{f}_1\mathbf{1}^T+\tau_{H_1}\mathbf{1}\mathbf{f}_2^T \\\hline
			\tau_{H_2}\mathbf{f}^T_1 & 0 & \tau_{H_1}\mathbf{f}^T_2\\\hline
			\tau_{H_1}\mathbf{f}_2\mathbf{1}^T+\tau_{H_2}\mathbf{1}\mathbf{f}_1^T & \tau_{H_1}\mathbf{f}_2  & \tau_{H_1}F_{\widetilde{H}_2}
		\end{array}\right],
	\end{align*}
	where $\mathbf{f}_1$ and $\mathbf{f}_2$ are the column vectors obtained from $\mathbf{f}_{H_1}^v$ and $\mathbf{f}_{H_2}^v$ by deleting the $v^\text{th}$ component (which is $0$), respectively. This implies that
	\begin{align*}
		&\mathbf{d}_G^TF_G\mathbf{d}_G\\
		=&\tau_{H_2}\mathbf{d}_{H_1}^TF_{H_1}\mathbf{d}_{H_1}+\tau_{H_1}\mathbf{d}_{H_2}^TF_{H_2}\mathbf{d}_{H_2}+4\tau_{H_2}m_{H_2}\mathbf{d}_{H_1}^T\mathbf{f}_{H_1}^v+4\tau_{H_1}m_{H_1}\mathbf{d}_{H_2}^T\mathbf{f}_{H_2}^v.
	\end{align*}
\end{proposition}


A formula for Kemeny's constant of a graph $G$ from \cite{kirkland2016kemeny} is 
\begin{equation}\label{formula:kemeny}
	\kappa(G)=\frac{\mathbf{d}_G^T F_G\mathbf{d}_G}{4m_G\tau_G}.
\end{equation}

In describing some graph families it will occasionally be convenient to have the following notation.
\begin{definition}
	Let $G_1, G_2$ be simple connected graphs and with labelled vertices $v_1 \in V(G_1)$, and
	$v_2 \in V(G_2)$. The \emph{1-sum} $G = G_1\bigoplus_{v_1,v_2} G_2$ is the graph created by 
	taking a copy of $G_1, G_2$, removing $v_1$, and replacing every edge of the form $\{i, v_1\} \in E(G_1)$
	with $\{i, v_2\}$. We often omit the subscript when the choice and/or labelling of vertices
	is clear. We say $G_1\bigoplus_vG_2$ has a \emph{1-separation}, and that $v$ is a \emph{1-separator} or \emph{cut vertex.}
\end{definition}
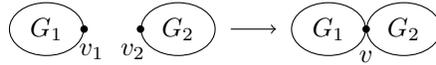
\begin{figure}[h]
	\centering
	\begin{tikzpicture}
		\tikzstyle{every node}=[circle, draw=none, fill=white, minimum width = 6pt, inner sep=1pt]
		\draw[] (1.75, 1)ellipse(14pt and 10pt);
		\draw[] (3.5, 1)ellipse(14pt and 10pt);
		\draw{
			(1.75,1)node[]{$G_1$}
			(3.5,1)node[]{$G_2$}
			(2.25,1)node[fill=black, minimum width = 2pt, label={[shift={(0.1,-.6)}]{$v_1$}}]{}
			(3,1)node[fill=black, minimum width = 2pt, label={[shift={(-0.1,-.6)}]{$v_2$}}]{}
		};
		
		\draw[->] (4.2,1) -- (4.8,1);
		
		\draw[] (5.5, 1)ellipse(14pt and 10pt);
		\draw[] (6.5, 1)ellipse(14pt and 10pt);
		\draw{
			(5.5,1)node[]{$G_1$}
			(6.5,1)node[]{$G_2$}
			(6,1)node[fill=black, minimum width = 2pt, label={[shift={(0,-.6)}]{$v$}}]{}
		};
	\end{tikzpicture}
	\caption{The graph $G=G_1\bigoplus_vG_2$ created from $G_1$ and $G_2$}
	\label{fig:vertexsum}
\end{figure}
For cleaner expressions we also introduce the following notation which was used in \cite{ciardo2022kemeny, faught20221}. 
\begin{definition}
	Let $G$ be a connected graph. The \emph{moment} of $v \in V$ is
	\begin{equation*}
		\mu(G, v) = \frac{1}{\tau_G} \mathbf{d}_G^T F_G \be_v.
	\end{equation*}
\end{definition}

\begin{example}\cite{kim2022families}\label{moment:P_n}
	Consider the path $\mathcal{P}_n=(1,2,\dots,n)$ where $n\geq 2$. Let $v$ be a vertex of $\mathcal{P}_n$. For $v=1,\dots,n$,
	\begin{align*}
		\kappa(\mathcal{P}_n)=\frac{1}{3}(n-1)^2+\frac{1}{6},\;\;\;\;\mu(\mathcal{P}_n,v)=(v-1)^2+(n-v)^2.
	\end{align*}
\end{example}

\begin{example}\cite{kim2022families}\label{moment:S_n}
	Consider a star $\mathcal{S}_n$ of order $n$ where $n\geq 3$. Suppose that $n$ is the center vertex. For $v=1,\dots,n$,
	\begin{align*}
		\kappa(\mathcal{S}_n) = n-\frac{3}{2},\;\;\;\;\mu(\mathcal{S}_n,v)=\begin{cases*}
			n-1, & \text{if $v=n$,}\\
			3n-5, & \text{if $v\neq n$.}
		\end{cases*}
	\end{align*}
\end{example}

\section{Braess edges on a path}\label{sec:path}

In this section, we consider which non-edges in a path are Braess edges. 
Throughout this section, we assume $k\geq 3$. 
Let $\cP_k=(1,\dots,k)$ be the path of length $k-1$. From Example~\ref{moment:P_n} with \eqref{formula:kemeny}, 
\begin{align}\label{eq:dfd(P)}
	\bd_{\cP_k}^TF_{\cP_k}\bd_{\cP_k}=\frac{4k^3-12k^2+14k-6}{3}=\frac{4}{3}(k-1)^3+\frac{2}{3}(k-1). 
\end{align} 

Throughout this paper, we will always assume that \( i<j \) when representing the edge $\{i,j\}$.
Let $\cP_k^{\{ i,j \}}$ be the graph obtained from $\cP_n$ by inserting the edge between vertex $i$ and vertex $j$ with $j-i\neq 1$. 
We now compute $\kappa(\cP_k^{\{ i,j \}})$ using Proposition~\ref{Prop:dFd with a cut-vertex}.
The vertex $i$ is a cut vertex in $\cP_k^{\{ i,j \}}$. 
There are exactly two branches of $\cP_k^{\{ i,j \}}$ at the vertex $i$:
one is the path $\cP_i$ of length $i-1$, and the other is the lollipop graph $\cL_{j-i+1,k-i+1}$,
where the \emph{lollipop graph} \( \cL_{a,b} \) is the graph obtained from a cycle of length \( a \) by appending a path of length \( b-a \). See Figure~\ref{Figure:P^i,j}.
\begin{figure}[h!]
	\begin{center}
		\begin{tikzpicture}[scale = .6]
			\tikzset{enclosed/.style={draw, circle, inner sep=0pt, minimum size=.10cm, fill=black}}
			
			\node[enclosed, label={below: $1$}] (v_1) at (0,0) {};
			\node[enclosed, label={below: $2$}] (v_2) at (1,0) {};
			\node[ label={below: $\cdots$}] (v_3) at (2,0) {};
			\node[enclosed, label={below: $i$}] (v_4) at (3,0) {};
			\node[ ] (v_5) at (4,0) {};
			\node[ label={below: $\cdots$}] (v_6) at (5,0) {};
			\node[] (v_7) at (6,0) {};
			\node[enclosed, label={below: $j$}] (v_8) at (7,0) {};
			\node[label={below: $\dots$}] (v_9) at (8,0) {};
			\node[enclosed, label={below: $k$}] (v_10) at (9,0) {};
			
			\draw [bend left = 40] (v_4) to (v_8);
			\draw (v_1) -- (v_10);
			\node at (4.5,-1.5) {\( \cP_k^{\{ i,j \}} \)};
			
		\end{tikzpicture}
		
		\begin{tikzpicture}[scale = .6]
			\tikzset{enclosed/.style={draw, circle, inner sep=0pt, minimum size=.10cm, fill=black}}
			
			\node[enclosed, label={below: $1$}] (v_0) at (-1,0) {};
			\node[enclosed, label={below: $2$}] (v_1) at (0,0) {};
			\node[ label={below: $\cdots$}] (v_2) at (1,0) {};
			\node[enclosed, label={below: $i$}] (v_3) at (2,0) {};
			\node[enclosed, label={below: $i$}] (v_4) at (3,0) {};
			\node[ ] (v_5) at (4,0) {};
			\node[ label={below: $\cdots$}] (v_6) at (5,0) {};
			\node[] (v_7) at (6,0) {};
			\node[enclosed, label={below: $j$}] (v_8) at (7,0) {};
			\node[label={below: $\dots$}] (v_9) at (8,0) {};
			\node[enclosed, label={below: $k$}] (v_10) at (9,0) {};
			
			\draw [bend left = 40] (v_4) to (v_8);
			\draw (v_0) -- (v_3);
			\draw (v_4) -- (v_10);
			\node at (0.5,-1.5) {\( B_1 := \cP_{i} \)};
			\node at (6,-1.5) {\( B_2 := \cL_{j-i+1,k-i+1} \)};
			
		\end{tikzpicture}
	\end{center}
	\caption{An illustration for \( \cP_k^{\{ i,j \}}, B_1 \) and \( B_2 \).}\label{Figure:P^i,j}
\end{figure}
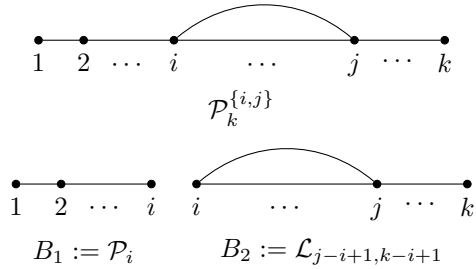

We let $s=i-1$, $l=j-i+1$, $B_1=\cP_{s+1}$, and $B_2=\cL_{l,k-s}$. 
The graph \( \cP_k^{\{ i,j \}} \) is obtained from \( B_1 \) and \( B_2 \) by identifying the vertex \( i \) in \( B_1 \) and the vertex \( i \) in \( B_2 \).
We already have \( \bd_{B_1}^TF_{B_1}\bd_{B_1} \) from \eqref{eq:dfd(P)} with \( k = s+1 \).
To compute \( \kappa(\cP_k^{\{ i,j \}}) \), we need to find \( \bd_{B_1}^T\mathbf{f}_{B_1}^i,  \bd_{B_2}^TF_{B_2}\bd_{B_2}\), and \( \bd_{B_2}^T\mathbf{f}_{B_2}^i \).

From \cite{kim2022families}, we have $\bd_{B_1}^T\mathbf{f}_{B_1}^i=s^2$. 
It appears in \cite{kirkland2016kemeny} that 
$$
\bd_{B_2}^TF_{B_2}\bd_{B_2}=\frac{2l(2(k-s)^3-4(k-s)l^2+3l^3-(k-s))}{3}.
$$
Furthermore, it can be seen that 
\begin{align*}
	\bd_{B_2}^T\mathbf{f}_{B_2}^i&=\frac{1}{3}(l-1)l(l+1)+(k-s-l)(l(k-s-l)+2l-2)\\
	&=l(k-s)^2-2(l^2-l+1)(k-s)+\frac{1}{3}(4l^3-6l^2+5l).
\end{align*}
Therefore, by Proposition \ref{Prop:dFd with a cut-vertex}, we have
\begin{align*}
	&\bd_{\cP_k^{\{ i,j \}}}^TF_{\cP_k^{\{ i,j \}}}\bd_{\cP_k^{\{ i,j \}}}\\
	=&~\tau_{B_2}\bd_{B_1}^TF_{B_1}\bd_{B_1}+\tau_{B_1}\bd_{B_2}^TF_{B_2}\bd_{B_2}+4\tau_{B_2}m_{B_2}\bd_{B_1}^T\mathbf{f}_{B_1}^i+4\tau_{B_1}m_{B_1}\bd_{B_2}^T\mathbf{f}_{B_2}^i\\
	=&~\frac{2}{3}\left( 12(l^2-l+1)s^2 + 12(l^3 - l^2k - l^2 + lk + l - k)s + l(3l^3 -4l^2k +2k^3 -k)\right).
\end{align*}
Hence, we have
\begin{align}\label{eq:k(P^ij)}
	\kappa\left(\cP_k^{\{ i,j \}}\right) 
	&= \frac{\bd_{\cP_k^{\{ i,j \}}}^TF_{\cP_k^{\{ i,j \}}}\bd_{\cP_k^{\{ i,j \}}}}{4kl} \\
	\notag &= \frac{1}{6kl}( 12(l^2-l+1)s^2 + 12(l^3 - l^2k - l^2 + lk + l - k)s \\
	\notag &+ l(3l^3 -4l^2k +2k^3 -k)). 
\end{align}
By computation, it is straightforward from Example~\ref{moment:P_n} and \eqref{eq:k(P^ij)} to obtain
\begin{align*}
	\kappa\left(\cP_k^{\{ i,j \}}\right) - \kappa\left(\cP_k\right) = \frac{f(s,l,k)}{6lk},
\end{align*}
where
\begin{align}\label{eq:f(s,l,k)}
	f(s,l,k) =  & ~4lk^2 -4( l + 3s - 3ls + 3l^2s + l^3)k \\
	\notag &+ s(12l^3 - 12l^2 + 12l) + s^2(12l^2 - 12l + 12) + 3l^4. 
\end{align}
This means that the non-edge \( ((s+1), (s +l)) \) is a Braess edge on \( \cP_k \) if 
\( f(s,l,k) >0 \). 

We now count the number of Braess edges on the path \( \cP_k \).

\begin{theorem}\label{thm:B(P)}
	For an integer \( k \ge 3 \), the number of Braess edges on the path \( \cP_k \) is
	\begin{align}\label{eq:Braess edges on path2}
		2\sum_{l= 3}^{\lfloor \sqrt{k} \rfloor} \left(\left\lfloor A(l) \right\rfloor +1 \right) + 2\max\left\{\left\lfloor A(\lfloor \sqrt{k} \rfloor+1) \right\rfloor+1 ,0\right\},
	\end{align}
	where 
	\begin{align*}
		A(l) := \frac{k-l}{2} - \sqrt{\frac{(3l^2-7l+3)k^2-2(l^3-3l^2+l)k-3l^2(l-1)}{12(l^2-l+1)}}.
	\end{align*}
\end{theorem}
\begin{proof}
	We rewrite the function \( f(s,l,k) \) in \eqref{eq:f(s,l,k)} as
	\begin{align*}
		12(l^2-l+1)\left(s-\frac{k-l}{2}\right)^2 -(3l^2-7l+3)k^2+2(l^3-3l^2+l)k+3l^2(l-1).
	\end{align*}
	Recall that a pair \( (s,l) \) gives a Braess edge on \( \cP_k \) if \( f(s,l,k)>0 \), equivalently,
	\begin{align}\label{eq:symmetry on P_k}
		s > (k-l) - A(l)\quad \mbox{ or}\quad s < A(l),
	\end{align}
	where
	\begin{align*}
		A(l) = \frac{k-l}{2} - \sqrt{\frac{(3l^2-7l+3)k^2-2(l^3-3l^2+l)k-3l^2(l-1)}{12(l^2-l+1)}}.
	\end{align*}
	By the symmetry in 	\eqref{eq:symmetry on P_k}, a pair \( (s,l) \) gives a Braess edge on \( \cP_k \) if and only if the pair \( \left((k-l)-s,l\right) \) gives a Braess edge on \( \cP_k \). 
	Therefore, the number of Braess edges for given \( l \) is 
	\begin{align*}
		2\times \max\left\{\left\lfloor A(l) \right\rfloor +1,0\right\}.
	\end{align*}
	Hence, the number of Braess edges on \( \cP_k \) is
	\begin{align}\label{eq:Braess edges on path1}
		2\sum_{l\ge 3}\max\left\{\left\lfloor A(l) \right\rfloor +1,0\right\}.
	\end{align}
	
	We now find the range of \( l \) such that the function \( A(l) \) is nonnegative, for given \( k \).
	We have
	\begin{align*}
		A(l) \ge 0 &\Longleftrightarrow
		\frac{(3l^3 -4kl^2 +4k^2-4k)l}{12(l^2-l+1)}\ge 0\\
		&\Longleftrightarrow 3l^3 -4kl^2 +4k^2-4k\ge 0.
	\end{align*}
	Let \( h(l)= 3l^3 -4kl^2 +4k^2-4k \).
	For \( k\ge 3 \), one can check that  \( h(\sqrt{k}) >0 \), \( h(\sqrt{k}+1)<0 \). Since \( h'(l) = 9l^2 -8kl \),
	we observe that \( h(l) \) is decreasing for \( 0\le l \le \sqrt{k} \).
	This ensures that \( A(l) \) is positive when \( 0< l \le \sqrt{k} \).
	Furthermore, we observe that \( A(l) \) is negative when \( \sqrt{k}+1 \le l \le k \), since \(h(l)\) is decreasing for \(\sqrt{k}+1 \le l \le 8k/9 \), and increasing for \(8k/9 \le l \le k\), in addition to \(h(\sqrt{k}+1)<0\) and \(h(k) = -k(k-2)^2<0 \).
	Therefore, we can express \eqref{eq:Braess edges on path1} as a finite sum, which completes the proof.
\end{proof}

We now show that if a non-edge crosses or contains a \emph{center vertex} $\lceil k/2\rceil$ in a path \( \cP_k \), then it is not a Braess edge.

\begin{theorem}\label{thm:nonedge crossing the center}
	For a path \( \cP_k \), we have \( \kappa\left(\cP_k^{\{ i,j \}}\right) - \kappa\left(\cP_k\right) < 0 \) if \( i\le\lceil k/2\rceil < j \).
\end{theorem}
\begin{proof}
	We know that \( \kappa\left(\cP_k^{\{ i,j \}}\right) - \kappa\left(\cP_k\right) = \frac{f(s,l,k)}{6lk}, \) where \( f(s,l,k) \) is defined in \eqref{eq:f(s,l,k)}, $s=i-1$, and $l=j-i+1$. To see \( f(s,l,k)<0 \) under the assumption, we rewrite \( f(s,l,k) \) as 
	\begin{align*} 
		f(s,l,k)=12(l^2-l+1)s^2+12(l^2-l+1)(l-k)s+l\cdot c(l,k), 
	\end{align*}
	where \( c(l,k)=4k^2-4(l^2+1)k+3l^3 \).
	We consider \( f(s,l,k) \) as a quadratic polynomial in variable \( s \), say \( f(s) \).  In this point of view, we give an equivalent condition of \( i\le\lceil k/2\rceil < j \) as follows: 
	\begin{align*} 
		0\le s\le \frac{k-1}{2}, \quad \frac{k}{2}-l\le s\le k-l,\text{ and } 3\le l\le k.
	\end{align*}
	\begin{figure}[h!]
		\centering
		\begin{tikzpicture}[domain=0.19:3.81]
			\draw[color=black] plot (\x,{(\x)^2 - 4*\x+4.3}) node[right] {$f(s)$} ;
			\draw (0,0) -- (4,0) node[right] {$s$};
			\draw[dotted] (0.2,0)  -- (0.2,0.2^2-4*0.2+4.3)  ;
			\draw[dotted] (3.8,0)  -- (3.8,3.8^2-4*3.8+4.3)  ;
			\draw[dotted] (1,0)  -- (1,1-4+4.3)  ;
			\draw[dotted] (2.9,0)  -- (2.9,2.9^2-4*2.9+4.3)  ;
			\draw (1,0.1) -- (1,-0.1) node[below]{ $\frac{k}{2}-l$}  ;
			\draw (2.9,0.1) -- (2.9,-0.1) node[below]{ $\frac{k-1}{2}$}  ;
			\draw (0.2,0.1) -- (0.2,-0.1) node[below]{ $0$} ;
			\draw (2,0.1) -- (2,-0.1) node[below]{ $\frac{k-l}{2}$};
			\draw (3.8,0.1) -- (3.8,-0.1) node[below]{ $k-l$} ;
			\draw[<->] (0.2,0.4) -- (1,0.4) ;
			\draw[<->] (2.9,0.4) -- (3.8,0.4) ;
			\draw (0.6,0.4) node[above]{$(*)$} ;
			\draw (3.35,0.4) node[above]{$(**)$} ;        
		\end{tikzpicture}
		\caption{A visualization of the first case with $k/2\ge l$. We show $(*)<(**)$ and hence the farthest point from $(k-l)/2$ is $s=k/2-l$.}
		\label{fig:visualization of f(s)}
	\end{figure}
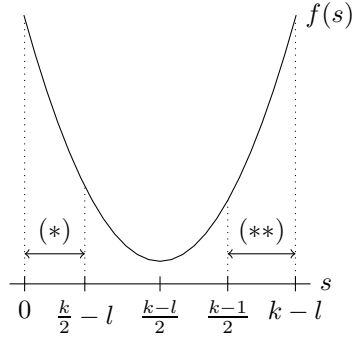
	
	The polynomial \( f(s) \) is minimized at \( s=\frac{k-l}{2}. \) Thus, \( f(s) \) is maximized when \( s \) is a point farthest from \(\frac{k-l}{2} \), which depends on \( k \) and \( l \). We consider two cases: (1) $\frac{k- 1}{2}\leq k-l$ and (2) $\frac{k- 1}{2}> k-l$. We claim that for both cases, the maximized value of \( f(s) \) is negative.
	\begin{enumerate}[leftmargin=*]
		\item  Suppose \( \frac{k-1}{2}\le k-l \), which is equivalent to \( 2l-1\le k \). First, let $k/2\ge l$. Then we have \( k-l-(k-1)/2-(k/2-l)=1/2 >0 \). In this case, the farthest point from \( (k-l)/2 \) is \( s=k/2-l \), see Figure~\ref{fig:visualization of f(s)}. For a given \( l \), let 
		\begin{align*}
			\widetilde{f}(k)&=f\left(\frac{k}{2}-l,l,k\right)\\
			&=-(3l^2-7l+3)k^2+2l(l^2-3l+1)k+3l^4. 
		\end{align*}
		For $l\ge3$, $\widetilde{f}(k)$ is decreasing for $k>l(l^2-3l+1)/(3l^2-7l+3)$. Since $k\ge 2l > l(l^2-3l+1)/(3l^2-7l+3)$, the possible maximum is attained at $k=2l$ and $\widetilde{f}(2l)=-l^2(5l^2-16l+8)<0$ for $l\ge3$. Second, let $k/2<l$. From the assumption, we have $k=2l-1$. Obviously, the farthest point from $(k-l)/2$ is $s=0$. For a given $l\ge3$, $f(0,l,2l-1)=-l(l-2)(5l^2-10l+4)<0.$
		
		\item Suppose \( \frac{k-1}{2}>k-l \). In this case, the farthest point from \( \frac{k-l}{2} \) is \( s=k-l \). For a given \( l \), let \( \hat{f}(k)=f(k-l,l,k)=4lk^2-4l(l^2+1)k+3l^4 \). For two roots \( r_1 \) and \( r_2 \) of \( \hat{f}(k) \) with \( r_1<r_2 \), we claim that \( r_1<l<2l-1<r_2 \) for \( l\ge3 \). When $l\ge3$, using the fact that $l^4-3l^3+2l^2+1>(l-1)^4$, we have
		\begin{align*}
			l-r_1=&~\frac{\sqrt{l^4-3l^3+2l^2+1}}{2}-\frac{(l-1)^2}{2}>0,  \\
			r_2-(2l-1)=&~\frac{\sqrt{l^4-3l^3+2l^2+1}}{2}+\frac{l^2-4l+3}{2}\\
			>&\frac{(l-1)^2}{2}+\frac{l^2-4l+3}{2}=(l-1)(l-2)>0.
		\end{align*}
		Additionally, since the coefficient of \( k^2 \) in \( \hat{f}(k) \) is positive,  we have \( \hat{f}(k)<0 \) for \( l\le k<2l-1 \). This completes the proof.  \qedhere
	\end{enumerate}
\end{proof}

\begin{remark}\label{rmk:far from the center}
	As we have seen in the proof of Theorem~\ref{thm:nonedge crossing the center}, \( f(s,l,k) \) is minimized at \( s=\frac{k-l}{2} \) for a given \( k \) and \( l \) and it increases as \( s \) moves away from \( \frac{k-l}{2} \). Roughly speaking, the further a non-edge moves away from the center of the path, the more likely it is to be a Braess edge.
\end{remark}

\begin{theorem}\label{thm:equiv Braess on path} 
	For a path \( \cP_k \), there exists a Braess edge making a cycle of length \( l\ge3 \) if and only if \( k>l^2-\frac{3}{4}l+\frac{7}{16}. \) 
\end{theorem}
\begin{proof}
	Let \( f(s,l,k) \) be defined as in \eqref{eq:f(s,l,k)}. By the symmetry of \( f(s,l,k) \) with respect to \( s=\frac{k-l}{2} \) together with Remark~\ref{rmk:far from the center}, there exists a Braess edge if and only if \( f(0,l,k)>0 \). The inequality \( f(0,l,k)=l\cdot (4k^2-4(l^2+1)k+3l^3)>0 \) is equivalent to
	\begin{align*} 
		r(l):=\frac{l^2+1}{2}+\frac{\sqrt{l^4-3l^3+2l^2+1}}{2}<k
	\end{align*}
	because we have $\frac{l^2+1}{2}-\frac{\sqrt{l^4-3l^3+2l^2+1}}{2}<l\le k$. One can check that \( l^2-\frac{3}{4}l+\frac{7}{16}>r(l) \) for \( l\ge3 \), which gives one direction of the proof. Now assume to the contrary that \( k\le l^2-\frac{3}{4}l+\frac{7}{16} \). Let \( h(l)=l^2-\frac{3}{4}l+\frac{7}{16}-r(l) \). One can check the following:
	\begin{align*} 
		h(3)<h(5)<h(4)<0.01308, \frac{dh(l)}{dl}<0 \text{ for } l\ge5, \text{ and } \lim_{l\rightarrow\infty}h(l)=0.
	\end{align*}
	Thus, \( h(l) \) is maximized at \( l=4 \) and \( h(l)<0.01308 \). Furthermore, the decimal part of $l^2-3l/4+7/16$ varies with $l$ and is presented in four different cases. Among these possibilities, the minimum value is greater than 0.18 and hence $\lfloor l^2- \frac{3}{4}l+\frac{7}{16}\rfloor=\lfloor l^2-\frac{3}{4}l+\frac{7}{16}-0.01308\rfloor$. Therefore, we have \[k\le \left\lfloor l^2-\frac{3}{4}l+\frac{7}{16}\right\rfloor = \left\lfloor l^2-\frac{3}{4}l+\frac{7}{16}-0.01308\right\rfloor \le \lfloor r(l) \rfloor \le r(l). \qedhere \]
\end{proof}

\subsection{Asymptotic behaviour}
We now examine asymptotic behavior of \( \sum_{l= 3}^{\lfloor \sqrt{k} \rfloor}\left\lfloor A(l) \right\rfloor \) described in Theorem~\ref{thm:B(P)}.

For $k>l>0$, we have 
\[   A(l) := A(l,k) = \frac{k-l}{2}- \sqrt{ \frac{(3l^2-7l+3)k^2 -2(l^3-3l^2+l)k-3l^2(l-1)}{12(l^2-l+1)} }.
\]
Also, let 
\[    A^*(k)   =  \sum_{l=3}^{\lfloor \sqrt{k}\rfloor}   \lfloor A(l)\rfloor .
\]

For convenience, we also introduce the notation
\begin{align*}
	h(l) = &~  \frac{3l^2-7l+3}{3(l^2-l+1)}  =   1-   \frac{4l}{3(l^2-l+1)} ,
	\\
	g(l)  = &~    \frac{2l(l^2-3l+1)}{3(l^2-l+1)} ,
	\\
	f(l) = &~ \frac{l^2(l-1)}{l^2-l+1} ,
\end{align*}   
so that we have
\begin{align}
	\label{eq.Aelldef}
	A(l)  =   \frac{k-l}{2}  -  \sqrt{ \frac{h(l)k^2}{4} - \frac{g(l)k}{4} -  \frac{f(l)}{4} }.
\end{align}

Recall that Euler's constant satisfies
\begin{equation}
	\label{eq.euler}
	\gamma   =  \lim_{N\rightarrow\infty}\left(-\ln N  +\sum_{l=1}^N\frac{1}{l}   \right)  
	=   0.5772\ldots
\end{equation}

\begin{proposition}
	\label{prop-asymp}
	\begin{equation}
		\label{eq.Alimprop}
		\lim_{k\rightarrow\infty}  \left(\frac{A^*(k)}{k} -\frac{1}{6}\ln(k)  \right)   =
		\frac{1}{2}S_1^{\infty}  +  \frac{1}{3}\gamma- \frac{2}{3}
	\end{equation}
	where $\gamma$ is Euler's constant and $S_1^{\infty}$ is the absolutely convergent sum
	\begin{equation}
		\label{eq.s1infdef}
		S_1^{\infty}=    
		\sum_{l=3}^{\infty}  \left(1- \sqrt{h(l)}  -  \frac{2}{3l} \right)  .
	\end{equation}
\end{proposition}

\begin{corollary}
	\label{cor-asymp}
	The number of Braess edges on the path $\cP_k$ is
	$\frac{k}{3}\ln k +\left(S_1^{\infty}+\frac{2}{3}\gamma-\frac{4}{3}\right)k + o(k)$.
\end{corollary}

\begin{remark}
	The value of $S_1^{\infty}$ is approximately $0.400$.  Thus the right-hand side of Equation (\ref{eq.Alimprop}) is approximately  $-0.274$.
\end{remark}

To prove Proposition \ref{prop-asymp}, we will need the following lemma.

\begin{lemma}
	\label{lem.taylor}
	Assume $D,D+u\geq w >0$.  Then 
	\begin{equation*}
		\label{eq.taylor}
		\left|  \sqrt{D+u}- \left( \sqrt{D}+ \frac{u}{2\sqrt{D}}\right) \right|  \leq   \frac{u^2}{8 w^{3/2}} .
	\end{equation*}
\end{lemma}
\noindent
\textbf{Proof of Lemma \ref{lem.taylor}:}  Let $\phi(x) =x^{1/2}$.  Then 
\[   \phi'(x) =\frac{1}{2} x^{-1/2}  \hspace{5mm}\hbox{and}\hspace{5mm}
\phi''(x) = -\frac{1}{4} x^{-3/2} .
\]
By Taylor's Theorem, there exists a $v$ between $D$ and $D+u$ such that 
\[    \phi(D+u)  =   \phi(D)  +  u\phi'(D)  +  \frac{u^2 \phi''(v)}{2}  .
\]
We have $|\phi''(v)|\leq 1/(4w^{3/2})$ because $v\geq w$, and the lemma follows.
\hfill $\Box$

\noindent
\textbf{Proof of Proposition \ref{prop-asymp}:}
First observe from Equation (\ref{eq.Aelldef}) that 
\begin{align*}
	\label{eq.Aelldef2}
	\frac{A(l)}{k} =  
	\frac{1}{2}\left(  1 - \frac{l}{k} -  \sqrt{h(l) - \frac{g(l)}{k}  - \frac{f(l)}{k^2} 
	}   \right) .
\end{align*}
Then we have
\begin{equation}
	\label{eq.bigsum}
	\frac{A^*(k)}{k}-\frac{1}{6}\ln k =\frac{1}{2}\left(S_0(k)+S_1(k) -S_2(k)+S_3(k) 
	+E(k) \right)
\end{equation}
where
\begin{align}
	\nonumber
	S_0(k) = & ~
	\frac{1}{k}\sum_{l=3}^{\lfloor \sqrt{k}\rfloor} \left( \lfloor A(l)\rfloor - A(l) \right)  ,
	\\
	\nonumber
	S_1(k)  = & ~
	\sum_{l=3}^{\lfloor \sqrt{k}\rfloor}   \left( 1-\sqrt{h(l)} - \frac{2}{3l}   \right) ,
	\\
	\nonumber
	S_2(k) = & ~
	\sum_{l=3}^{\lfloor \sqrt{k}\rfloor}  \frac{l}{k} ,  
	\\
	\label{eq.S3def}
	S_3(k) = & ~
	\sum_{l=3}^{\lfloor \sqrt{k}\rfloor}   
	\left( \sqrt{h(l)} - \sqrt{ h(l)- \frac{g(l)}{k} - \frac{f(l)}{k^2} } \right) ,  \quad  \hbox{and}
	\\
	\label{eq.Edef}
	E(k) = &  ~ \frac{2}{3}\left( -\ln \sqrt{k} + \sum_{l=3}^{\lfloor \sqrt{k}\rfloor} \frac{1}{l} \right).
\end{align}   

Since $-1\leq \lfloor A(l)\rfloor - A(l)\leq 0$, we see that $-\sqrt{k}/k\leq S_0(k)\leq 0$, and 
hence 
\begin{equation}
	\label{eq.S0lim}
	\lim_{k\rightarrow\infty}S_0(k)  = 0 .
\end{equation}
The sum $S_2(k)$ is easy to evaluate exactly, and we conclude
\begin{equation}
	\label{eq.S2lim}
	\lim_{k\rightarrow\infty}S_2(k)  = \frac{1}{2} .
\end{equation}
For $E(k)$, notice that the lower limits of the sums are different in 
Equations (\ref{eq.euler}) and (\ref{eq.Edef}).  Accounting for this, and using the fact that 
$\lim_{k\rightarrow \infty}(\ln\sqrt{k}-\ln \lfloor \sqrt{k}\rfloor)=0$, we see that
\[   
\lim_{k\rightarrow\infty}E(k)+\frac{2}{3}\sum_{l=1}^2\frac{1}{l}  = \frac{2}{3}\gamma ,  
\]
and hence
\begin{equation}
	\label{eq.Elim}
	\lim_{k\rightarrow\infty}E(k)  = \frac{2}{3}\gamma - 1.
\end{equation}

At this point, it is useful to have some simple bounds on the functions $h$, $g$, and $f$.
For $l\geq 3$, we have 
\begin{align}
	\label{eq.hbd1a}
	h(l)  = &~ 1 - \frac{4l}{3(l^2-l+1)}  >   1- \frac{4l}{3(l^2-l)}
	=   1- \frac{4}{3(l-1)}  
	\\
	\label{eq.hbd1b}
	\geq &~  1-\frac{4}{6} = \frac{1}{3}  ,
	\\
	\nonumber
	g(l)  = &~  \frac{2l(l^2-3l+1)}{3(l^2-l+1)}       <  \frac{2l}{3} ,
	\hspace{5mm}\hbox{and}
	\\
	\label{eq.fbd1}
	f(l)  = & ~  \frac{l^2(l-1)}{l^2-l+1}     <   \frac{l^2(l-1)}{l^2-l}  
	=  l.
\end{align}
Also, 
\begin{align}
	\label{eq.hbd2}
	\frac{1-h(l)}{2}-\frac{2}{3l} =  \frac{2l}{3(l^2-l+1)}-\frac{2}{3l} 
	=  \frac{2(l-1)}{3l(l^2-l+1)}  
\end{align}
and hence
\begin{align}
	\label{eq.hbd3}
	0  \leq   \frac{1-h(l)}{2}-\frac{2}{3l} \leq   \frac{2(l-1)}{3l(l^2-l)}  =  \frac{2}{3l^2} .  
\end{align}

To handle $S_1(k)$, we shall apply Lemma \ref{lem.taylor} with $D=1$, $u=h(l)-1$, and $w=1/4$.
The bound (\ref{eq.hbd1b}) shows that $D+u>w$, so we have
\begin{align*}
	\left| 1-\sqrt{h(l)} - \frac{2}{3l}   \right|   = &~ 
	\left|\sqrt{h(l)} -1  -\frac{h(l)-1}{2}  +  \frac{h(l)-1}{2}  + \frac{2}{3l}   \right|
	\\
	\leq  & ~
	\left|\sqrt{h(l)} -1  -\frac{h(l)-1}{2}\right|  +  \left|\frac{1-h(l)}{2}  - \frac{2}{3l}   \right| 
	\\
	\leq &  ~ \frac{ (h(l)-1)^2}{8(1/4)^{3/2}}  +  \frac{2}{3l^2} 
	\\
	&~ \hspace{18mm}\hbox{(by Lemma \ref{lem.taylor} and Equation (\ref{eq.hbd3}))}
	\\
	\leq & ~\frac{16}{9(l-1)^2}   +  \frac{2}{3l^2} 
	\hspace{9mm}\hbox{(by Equation (\ref{eq.hbd1a}))}.        
\end{align*}        
It follows that  the series $S_1^{\infty}$ of Equation (\ref{eq.s1infdef}) is absolutely convergent, 
and of course
\begin{align}
	\label{eq.S1lim}
	\lim_{k\rightarrow\infty}S_1(k)  = S_1^{\infty} .
\end{align}      

It remains to deal with $S_3(k)$.    We begin by applying Lemma \ref{lem.taylor} with 
\[     D=h(l)   \hspace{5mm}\hbox{and}\hspace{5mm}  u =  - \frac{g(l)}{k} - \frac{f(l)}{k^2} .
\]   
By Equations (\ref{eq.hbd1a}--\ref{eq.fbd1}), when $3\leq l \leq \sqrt{k}$, we have $D>\frac{1}{3}$ and 
\begin{align}
	\label{eq.hgfbd2}
	0  >    u  >  -\frac{2l}{3k}  - \frac{l}{k^2}  \geq  -\frac{2}{3\sqrt{k}}
	- \frac{1}{k^{3/2}}   \geq   - \frac{5}{3\sqrt{k}}.
\end{align}
Since we are really interested in $k$ tending to infinity, it is harmless to assume that $k\geq 36$.  Then
Equation (\ref{eq.hgfbd2}) implies that 
\begin{align*}
	\label{eq.hgfbd3}
	D+u  > \frac{1}{3} - \frac{5}{3\sqrt{36}}    = \frac{1}{18} 
	\hspace{5mm}\hbox{when $3\leq l \leq \sqrt{k}$ and $k\geq 36$.}
\end{align*}
Therefore we can take $w=1/18$ in Lemma \ref{lem.taylor}.
Thus the lemma tells us that we can write
\begin{equation}
	\label{eq.hdiff1}
	\sqrt{h(l)} - \sqrt{ h(l)- \frac{g(l)}{k} - \frac{f(l)}{k^2} } =  
	\frac{  \frac{g(l)}{k}+\frac{f(l)}{k^2}  }{2\sqrt{h(l)} } +\epsilon(l,k).
\end{equation}
where the remainder term $\epsilon$ satisfies
\begin{align}
	\label{eq.epsbound}
	\left|\epsilon(l,k) \right|   \leq      \frac{\left(\frac{5}{3\sqrt{k}}\right)^2}{8\left(\frac{1}{18}\right)^{3/2}}
	=  \frac{75}{2\sqrt{2}\,k} 
	\hspace{5mm}\hbox{when $3\leq l \leq \sqrt{k}$ and $k\geq 36$.}               
\end{align}

Now, by Equations (\ref{eq.S3def}) and (\ref{eq.hdiff1}), we can rewrite $S_3(k)$ as
\begin{align*}
	\label{eq.Sdecomp}
	S_3(k) = S_{3,g}(k) +S_{3,f}(k)+ S_{3,\epsilon}(k),
\end{align*}
where
\begin{align*}
	S_{3,g}(k)=~\frac{1}{k} \sum_{l=3}^{\lfloor \sqrt{k}\rfloor}\frac{g(l)}{2\sqrt{h(l)}},\;\;\; S_{3,f}(k) = \frac{1}{k^2} \sum_{l=3}^{\lfloor \sqrt{k}\rfloor}\frac{f(l)}{2\sqrt{h(l)}},\;\;\; S_{3,\epsilon}(k)  = \sum_{l=3}^{\lfloor \sqrt{k}\rfloor}  \epsilon(l,k).
\end{align*}

By Equation (\ref{eq.epsbound}), $|S_{3,\epsilon}(k)| \leq  \frac{75}{2\sqrt{2k}}$ when $k\geq 36$, and so
\begin{equation}
	\label{eq.s3epslim}
	\lim_{k\rightarrow\infty}  S_{3,\epsilon}(k)  = 0.
\end{equation}
By Equations (\ref{eq.hbd1b}) and (\ref{eq.fbd1}), we see that 
\[   |S_{3,f}(k)| \leq  \frac{1}{k^2} \frac{\sqrt{k}(\sqrt{k}+1)}{4\sqrt{1/3}} ,
\] 
and therefore
\begin{equation}
	\label{eq.s3flim}
	\lim_{k\rightarrow\infty}  S_{3,f}(k)  = 0.
\end{equation}
Next, we claim that 
\begin{equation}
	\label{eq.s3glim}
	\lim_{k\rightarrow\infty}  S_{3,g}(k)  = \frac{1}{6} . 
\end{equation}
To obtain this result, we write $g(l)$ as
\begin{align*}  
	g(l)  = \frac{2l(l^2-3l+1)}{3(l^2-l+1)} =  \frac{2(l-2)}{3} + g_1(l),\;\;\hbox{where}\;\;\;   g_1(l) = -\frac{4(l-1)}{3(l^2-l+1)}  .     
\end{align*}
Accordingly, we write
\begin{align}
	\label{eq.gg1suma}   
	S_{3,g}(k)   =  T(k) + T_1(k), 
\end{align}
where
\begin{align*}
	T(k) =  \frac{1}{k} \sum_{l=3}^{\lfloor \sqrt{k}\rfloor}  \frac{(l-2)}{3\sqrt{h(l)}}   
	\;\;\;\hbox{and}\;\;\;T_1(k)  = \frac{1}{k} \sum_{l=3}^{\lfloor \sqrt{k}\rfloor}  \frac{g_1(l)}{2\sqrt{h(l)} }.
\end{align*}
We have from Equation (\ref{eq.hbd1b}) that 
\begin{align*}
	\left| \frac{g_1(l)}{2\sqrt{h(l)}} \right|  \leq   \frac{2(l-1)}{3(l^2-l)\sqrt{1/3}} 
	<   \frac{4}{3l} ,
\end{align*}
and therefore  
$\left| T_1(k) \right|  \leq   \frac{4}{3}\,\sqrt{k}/k$,
from which we obtain
\begin{align}
	\label{eq.T1lim}
	\lim_{k\rightarrow\infty}  T_1(k)  = 0.
\end{align}
Now we consider $T(k)$.  Note that $h(l)<1$ for every $l\geq 3$.
Let $\delta\in (0,1)$.  Since $\lim_{l\rightarrow\infty}h(l)=1$, we can choose $L$ such that 
$1> \sqrt{h(l)} > 1-\delta$ whenever $l>L$.  Then
\begin{equation}
	\label{eq.Tbd1}
	T(k)   \geq   \frac{1}{k} \sum_{l=3}^{\lfloor \sqrt{k}\rfloor}  \frac{(l-2)}{3}   = 
	\frac{ ( \lfloor \sqrt{k}\rfloor-1)(\lfloor \sqrt{k}\rfloor-2)}{6k} .
\end{equation}
For the corresponding upper bound, assume $\sqrt{k}>L$.  Then
\begin{eqnarray}
	\nonumber
	T(k)  & \leq &   \frac{1}{k} \sum_{l=L}^{\lfloor \sqrt{k}\rfloor}  \frac{(l-2)}{3(1-\delta)} 
	+  \frac{1}{k}\sum_{l=3}^{L-1} \frac{l-2}{3\sqrt{h(l)}}
	\\
	& \leq & 
	\label{eq.Tbd2}
	\frac{ ( \lfloor \sqrt{k}\rfloor-1)(\lfloor \sqrt{k}\rfloor-2)}{6k(1-\delta)} 
	+  \frac{1}{k}\sum_{l=3}^{L-1} \frac{l-2}{3\sqrt{h(l)}}.
\end{eqnarray}
Combining Equations (\ref{eq.Tbd1}) and (\ref{eq.Tbd2}) yields (since $L$ is fixed)
\[   \frac{1}{6}  \leq   \liminf_{k\rightarrow\infty}T(k)  \leq   \limsup_{k\rightarrow\infty}T(k) 
\leq   \frac{1}{6(1-\delta)}.
\]
Since the above holds for every $\delta$ in $(0,1)$, we conclude
that
\begin{equation}   
	\label{eq.Tlim1}
	\lim_{k\rightarrow\infty}T(k) = \frac{1}{6} .                  
\end{equation}
Now the claimed Equation (\ref{eq.s3glim}) follows from Equations (\ref{eq.gg1suma}),  
(\ref{eq.T1lim}), and (\ref{eq.Tlim1}).  Then by Equations (\ref{eq.s3epslim})--(\ref{eq.s3glim}), we obtain
\begin{equation}   
	\label{eq.S3lim}
	\lim_{k\rightarrow\infty}S_3(k) = \frac{1}{6} .                  
\end{equation}

Finally, by Equations (\ref{eq.bigsum}), (\ref{eq.S0lim})--(\ref{eq.Elim}), (\ref{eq.S1lim}), and 
(\ref{eq.S3lim}), we conclude
\begin{eqnarray*}
	\lim_{k\rightarrow\infty}  \left(\frac{A^*(k)}{k} -\frac{1}{6}\ln(k)  \right)  
	& = & \frac{1}{2}\left(   0+S_1^{\infty} -\frac{1}{2} +\frac{1}{6} +\frac{2}{3}\gamma-1\right)
	\\
	& = & 
	\frac{1}{2}S_1^{\infty}  +  \frac{1}{3}\gamma- \frac{2}{3}.
\end{eqnarray*}
This proves the proposition.   
\hfill $\Box$

\section{Braess edges on a spider graph}\label{sec:spider}

First we mention a calculation that is useful for identifying Braess edges. The result in \cite[Theorem~4.2]{faught20221} can be recast as follows.
\begin{theorem}\label{thm:Braess with cut vertex}
	Let $G=G_1\bigoplus_v G_2$. Let $\tilde{G}$ (resp. $\tilde{G}_1$) be the graph obtained from $G$ (resp. $G_1$) by adding an edge to $G_1$. Then,
	\begin{align*}
		&\kappa(\tilde{G})-\kappa(G)\\
		=&\kappa(\tilde{G}_1)-\kappa(G_1)+\frac{m_{G_2}}{m_{\tilde{G}}}\left(\mu(\tilde{G}_1,v)-\kappa(\tilde{G}_1)-\mu(G_2,v)+\kappa(G_2)\right)\\
		&-\frac{m_{G_2}}{m_G}\left(\mu(G_1,v)-\kappa(G_1)-\mu(G_2,v)+\kappa(G_2)\right).
	\end{align*}
\end{theorem}

A \textit{spider graph} is a tree with one vertex of degree at least $3$ and all other vertices of degree $1$ or $2$. For $a\geq 1$ and $b\geq 2$, we use $\mathcal{S}_{a,b}$ to denote the spider graph with a vertex $v$ of degree $b$ such that each pendent vertex is at distance $a$ from $v$. From \cite[Proposition~6.3]{ciardo2022kemeny}, we can find
\begin{align*}
	\kappa(\mathcal{S}_{a,b}) = \left(b-\frac{2}{3}\right)a^2+\frac{1}{6}.
\end{align*}
Considering the distance matrix of $\mathcal{S}_{a,b}$, we can find 
\begin{align*}
	\mu(\mathcal{S}_{a,b},v) = ba^2.
\end{align*}
Then, 
\begin{align}\label{eqn:temp111}
	\mu(\mathcal{S}_{a,b},v)-\kappa(\mathcal{S}_{a,b}) =\frac{2}{3}a^2-\frac{1}{6}= \mu(\mathcal{P}_{a+1},w)-\kappa(\mathcal{P}_{a+1})
\end{align}
where $w$ is a pendent vertex of the path. 

We remark that by \cite[Theorem~2.2]{ciardo2022kemeny}, every non-edge in $S_{1,b}$ for $b\geq 2$ is a Braess edge, and so we shall consider $S_{a,b}$ for $a \geq 2$ throughout this section.

\begin{proposition}
	Let $a\geq 2$, $b\geq 2$, and $v$ be the center vertex of $\mathcal{S}_{a,b}$. Choose two pendent vertices in $\mathcal{S}_{a,b}$. We let $\mathcal{P}_{2a+1}=(1,\dots,2a+1)$ be the path from one pendent vertex to the other. For $1\leq i<j\leq 2a+1$ with $j-i\geq 2$, let $\mathcal{S}_{a,b}^{\{ i,j \}}$ denote the graph obtained from $\mathcal{S}_{a,b}$ by adding edge $\{ i,j \}$. Then, we have
	\begin{align}\nonumber
		&\kappa(\mathcal{S}_{a,b}^{\{ i,j \}})-\kappa(\mathcal{S}_{a,b}) \\\label{exp:spider}
		= &\kappa(\cP_{2a+1}^{\{ i,j \}}) - \kappa\left(\cP_{2a+1}\right)+\frac{a(b-2)}{ab+1}\left(\mu(\cP_{2a+1}^{\{ i,j \}},v) -\kappa(\cP_{2a+1}^{\{ i,j \}})-\frac{2}{3}a^2+\frac{1}{6}\right).
	\end{align}
\end{proposition}
\begin{proof}
	When $b=2$, it is trivial. Let $b= 3$.   Choosing $G_1=\mathcal{S}_{a,2}$, $\tilde{G}_1 = \mathcal{S}_{a,2}^{\{ i,j \}}$, and $G_2 = \mathcal{P}_{a+1}$, using Theorem~\ref{thm:Braess with cut vertex} with \eqref{eqn:temp111}, Equation \eqref{exp:spider} follows. Similarly, for $b \geq 4$, choosing $G_1=\mathcal{S}_{a,2}$, $\tilde{G}_1 = \mathcal{S}_{a,2}^{\{ i,j \}}$, and $G_2 = \mathcal{S}_{a,b-2}$, we can establish the result.
\end{proof}

To find the moment $\mu(\cP_{2a+1}^{\{ i,j \}},v)$, we calculate each entry in $F_{\cP_{2a+1}^{\{ i,j \}}}$ case by case:
\begin{lemma} \label{lem:f_pq}
	Let $\mathcal{P}_k=(1,\dots,k)$ be the path of length $k-1$. 
	Let $F_{\mathcal{P}_k^{\{i, j\}}}=(f_{p,q})_{1\le p,q \le k}$ and $l$ be the length of the cycle in $\mathcal{P}_k^{\{ i,j \}}$. 
	We have the following:
	\begin{itemize}
		\item if $p,q<i$ or $p,q>j$, then $f_{p,q} = l|p-q|$;
		\item if $p<i$ and $i\le q\le j$, then $f_{p,q} = l(i-p)+(q-i)(l-q+i)$;
		\item if $p<i$ and $q>j$, then $f_{p,q} = l(i-p)+l(q-j)+(l-1)$;
		\item if $i\le p \le j$ and $i\le q \le j$, then $f_{p,q} = |p-q|(l-|p-q|)$;
		\item  if $i\le p\le j$ and $q>j$, then $f_{p,q}=l(q-j)+(j-p)(l-j+p)$.
	\end{itemize}
	Other cases not stated are obtained by symmetry.
\end{lemma}
\begin{proof}
	Using a similar argument as in \cite[Example~3.4]{kirkland2016kemeny}, one can easily compute $f_{p,q}$, considering two cases: the vertex $p$ is on the cycle in $\mathcal{P}_k^{\{i, j\}}$ if $i \le p \le j$, and $p$ is not on the cycle if $p<i$ or $p>j$.
\end{proof}

Preserving the notation of the above lemma, by the definition of the moment, we have
\begin{align} \nonumber
	\mu(\cP_{2a+1}^{\{ i,j \}},v) 
	=&~ \frac{1}{\tau_{\cP_{2a+1}^{\{ i,j \}}}} \mathbf{d}_{\cP_{2a+1}^{\{ i,j \}}}^T F_{\cP_{2a+1}^{\{ i,j \}}} e_v \\\label{eq:moment Pij}
	=&~ \frac{1}{l}\left( f_{1,v} + 2\sum_{m=2}^{i-1} f_{m,v} + 2\sum_{m=i}^{j} f_{m,v} + 2\sum_{m=j+1}^{2a} f_{m,v} + f_{2a+1,v} \right).
\end{align}

Note that $v$ is the center of $\cP_{2a+1}$, that is, $v=a+1$.
Let $i = s+1$ and $j = s+l$, so that $l$ is the length of the cycle in $\cP_{2a+1}^{\{i,j\}}$.
Now we shall find the moment $\mu(\cP_{2a+1}^{\{ i,j \}},v)$ by considering two cases: 
\begin{enumerate}[label=(\roman*)]
	\item $v$ is contained in the cycle of $\cP_{2a+1}^{\{ i,j \}}$;
	\item $v$ is not contained in the cycle of $\cP_{2a+1}^{\{ i,j \}}$.
\end{enumerate}

For the case (i), we have $i\le v \le j$, equivalently, $s+1\leq a+1\leq s+l$.
Applying Lemma~\ref{lem:f_pq} and \eqref{eq:moment Pij}, we obtain
\begin{multline}\label{eq:moment1}
	\mu(\cP_{2a+1}^{\{ i,j \}},v) = -\frac{2(2a-2l+3)}{l}s^2+\frac{2(2a-2l+3)(2a-l+1)}{l}s\\
	+\frac{4l^2}{3}-2(3a+2)l+10a^2+14a+\frac{14}{3}-\frac{2(2a^3+5a^2+4a+1)}{l}.
\end{multline}

For the case (ii), considering the symmetry, we may only consider $i<j<v$, equivalently, $0\leq s$ and $s+l\leq a$.
Again, applying Lemma~\ref{lem:f_pq} and \eqref{eq:moment Pij}, we obtain
\begin{align}\label{eq:moment2}
	\mu(\cP_{2a+1}^{\{ i,j \}},v)= -\frac{2(l^2-l+1)}{l}s+2a^2+2a-\frac{2}{3}l^2+\frac{2}{3}.
\end{align}

We identify some non-Braess edges on $\mathcal{S}_{a,b}$.

\begin{lemma}\label{lem:ineq1}
	Let $a\geq 2$. Suppose that the cycle of  $\cP_{2a+1}^{\{ i,j \}}$ contains $v$ where $v = a+1$. Then, $$\mu(\cP_{2a+1}^{\{ i,j \}},v)- \kappa\left(\cP_{2a+1}\right)-\frac{2}{3}a^2+\frac{1}{6}\leq 0.$$
\end{lemma}
\begin{proof}
    Let $s=i-1$ and $l=j-i+1$.
	Using \eqref{eq:moment1}, we can find
	\begin{align*}
		H_1(s,l,a):=&~l\left(\mu(\cP_{2a+1}^{\{ i,j \}},v)- \kappa\left(\cP_{2a+1}\right)-\frac{2}{3}a^2+\frac{1}{6}\right)\\
		=&-2(2a - 2l + 3)s^2 + 2(2a - 2l + 3)(2a - l + 1)s\\
		& - 4a^3 + (8l - 10)a^2  -2(l-1)(3l-4)a +\frac{2}{3}(l-1)(2l^2 - 4l + 3).
	\end{align*}
	
	Let $a\geq l-1$. If we consider $H_1$ as a polynomial in variable $s$, then $H_1$ is concave down and the maximum is given by
	\begin{align*}
		H_1\left(\frac{2a-l+1}{2},l,a\right)= -(l-1)^2a + \frac{1}{3}l^3 - \frac{1}{2}l^2 + \frac{2}{3}l -  \frac{1}{2},
	\end{align*}
	which is a linear function in variable $a$ with a negative slope. Since $H_1\left(\frac{l-1}{2},l,l-1\right)=-\frac{1}{6}(l-1)(4l^2-11l+3)<0$ for $l\ge 3$, we have $H_1(s,l,a)<0$ for $a\geq l$.
	
	Suppose $\frac{l-1}{2}<a<l-1$. Then, the polynomial $H_1$ in variable $s$ is concave up. Examining the proof of Theorem~\ref{thm:nonedge crossing the center}, given $l$ and $a$, the maximum is attained at $s=0$. So,
	\begin{align*}
		H_1\left(0,l,a\right)= - 4a^3+(8l - 10)a^2  -2(l-1)(3l-4)a + \frac{2}{3}(l-1)(2l^2 - 4l + 3).
	\end{align*}
	Taking the derivative of $H_1(0,l,a)$ with respect to $a$ yields $ - 12a^2  + 4(4l - 5)a -2(l-1)(3l-4)$, which can be shown to be negative. Since $H_1\left(0,l,\frac{l}{2}\right)=-\frac{1}{6}(l-3)(l-1)(l+2)\leq 0$ for $l \geq 3$, we have $H_1(s,l,a)<0$ for $\frac{l-1}{2}<a<l-1$. 
	
	Finally, we suppose $l =2a+1$. Then $s = 0$ and we have
	$$H_1\left(0,l,\frac{l-1}{2}\right) = -\frac{1}{6}l(l-1)(l-5).$$
	Since $a\geq 2$, we have $l\geq 5$. It can be seen that $H_1\left(0,l,\frac{l-1}{2}\right)\leq 0$ for $l\geq 5$. Therefore, the conclusion follows. 
\end{proof}

\begin{theorem}\label{prop:nonBraess1}
	Let $a\geq 2$ and $b\geq 2$. Suppose that the cycle in $\mathcal{S}_{a,b}^{\{ i,j \}}$ contains $v$, where $v$ is the center vertex of $\mathcal{S}_{a,b}$. Then, $$\kappa(\mathcal{S}_{a,b}^{\{ i,j \}})-\kappa(\mathcal{S}_{a,b}) < 0. $$
\end{theorem}
\begin{proof}
	Consider \eqref{exp:spider}. From Theorem~\ref{thm:nonedge crossing the center}, we have $\kappa(\cP_{2a+1}^{\{ i,j \}}) - \kappa\left(\cP_{2a+1}\right)<0$. If $\mu(\cP_{2a+1}^{\{ i,j \}},v) -\kappa(\cP_{2a+1}^{\{ i,j \}})-\frac{2}{3}a^2+\frac{1}{6}\leq 0$, then $\kappa(\mathcal{S}_{a,b}^{\{ i,j \}})-\kappa(\mathcal{S}_{a,b}) < 0$. Suppose that $\mu(\cP_{2a+1}^{\{ i,j \}},v) -\kappa(\cP_{2a+1}^{\{ i,j \}})-\frac{2}{3}a^2+\frac{1}{6}>0$. Since $0 \leq \frac{a(b-2)}{ab+1}<1$, we have
	\begin{align}\label{tempeqn:11}
		\kappa(\mathcal{S}_{a,b}^{\{ i,j \}})-\kappa(\mathcal{S}_{a,b}) < \mu(\cP_{2a+1}^{\{ i,j \}},v)- \kappa\left(\cP_{2a+1}\right)-\frac{2}{3}a^2+\frac{1}{6}.
	\end{align}
	By Lemma~\ref{lem:ineq1}, it follows that $\kappa(\mathcal{S}_{a,b}^{\{ i,j \}})-\kappa(\mathcal{S}_{a,b})<0$. This proves the theorem.
\end{proof}

\begin{figure}
	\begin{center}
		\begin{tikzpicture}
			\tikzset{enclosed/.style={draw, circle, inner sep=0pt, minimum size=.10cm, fill=black}}
			\def \radius {1cm}
			
			\node[enclosed] (origin) at (360:0mm) {};
			\foreach \i [count=\ni from 0] in {6,5,4,3,2,k,k-1}{
				\node[enclosed] at ({210-\ni*40}:\radius) () {};
				\node[enclosed] at ({210-\ni*40}:2*\radius) (u\ni) {};
				\draw (origin)--(u\ni);
			}
			\draw[thick, loosely dotted] (240:\radius) arc[start angle=240, end angle=300, radius=\radius];
		\end{tikzpicture}
	\end{center}
	\caption{There is no Braess edge in $\mathcal{S}_{2,b}$ for $b\geq 2$.}\label{Figure:2-star}
\end{figure}
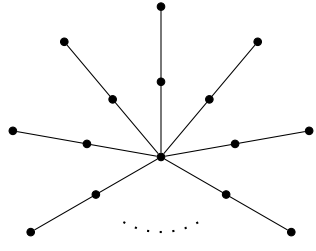

We obtain the following corollary for a family of trees with no Braess edges (see Figure~\ref{Figure:2-star}).

\begin{corollary}\label{cor:no Braess}
	For an integer $b\geq2$, $B(\mathcal{S}_{2,b}) = 0$. 
\end{corollary}

We now consider which edges are Braess on $\mathcal{S}_{a,b}$ in the following Lemma. 
\begin{lemma}\label{lem:Braess on spider}
	Let $a\geq 2$. If $\kappa(\cP_{2a+1}^{\{ i,j \}}) - \kappa\left(\cP_{2a+1}\right)>0$, then $$\mu(\cP_{2a+1}^{\{ i,j \}},v)-\kappa(\cP_{2a+1}^{\{ i,j \}})-\frac{2}{3}a^2+\frac{1}{6}>0.$$
	This implies that if $\{ i,j \}$ is a Braess edge for $\cP_{2a+1}$, then it also is for $\mathcal{S}_{a,b}$.
\end{lemma}
\begin{proof}
Let $s=i-1$ and $l=j-i+1$.	
 Using \eqref{eq:k(P^ij)} and \eqref{eq:moment2}, we can find that
	\begin{align*}
		F_1(s,l,a) :=&~ 6l(2a+1)\left(\mu(\cP_{2a+1}^{\{ i,j \}},v)-\kappa(\cP_{2a+1}^{\{ i,j \}})-\frac{2}{3}a^2+\frac{1}{6}\right)\\
		=& -12(l^2 - l + 1)s^2 - 12l(l^2 - l + 1)s + l(8a^2 + 12a - 3l^3 + 4).
	\end{align*}
	Since $\{ i,j \}$ is a Braess edge, we have from Theorem~\ref{thm:equiv Braess on path} and \eqref{eq:symmetry on P_k} that \begin{align*}
		a>K\;\;\text{and}\;\;\;\; 0\leq s< \frac{(2a+1)-l}{2}-\sqrt{\frac{N}{D}}:=s_0, 
	\end{align*}
	where  
	\begin{align*}
		K =&~\frac{1}{2}\left(l^2-\frac{3}{4}l-\frac{9}{16}\right),\\
		N =&~ (3l^2-7l+3)(2a+1)^2-2(l^3-3l^2+l)(2a+1)-3l^2(l-1), \\
		D =&~ 12(l^2-l+1).
	\end{align*}
	Considering $F_1(s,l,a)$ as a polynomial in variable $s$, observe that it is concave down. It is enough to show that $F_1(0,l,a)>0$ and $F_1(s_0,l,a)>0$ for $l\geq 3$. For the first inequality, let $R := R(l,a) = l(8a^2 + 12a - 3l^3 + 4)$. We can see that 
	\begin{align*}
		F_1(0,l,a)= R(l,a) > R(l,K) = 2l^5 - 6l^4 + \frac{39}{8}l^3 - \frac{45}{16}l^2 + \frac{161}{128}l>0.
	\end{align*}
	Now we consider
	\begin{align*}
		F_1(s_0,l,a) =&~(2a+1)\sqrt{ND}+R-D\left(\frac{(2a+1-l)^2}{4}+\frac{l(2a+1-l)}{2}+\frac{N}{D}\right).
	\end{align*}
	Then $F_1(s_0,l,a)>0$ is implied by $$F_2(l,a) := (2a+1)^2ND-\left(R-D\left(\frac{(2a+1-l)^2}{4}+\frac{l(2a+1-l)}{2}+\frac{N}{D}\right)\right)^2>0.$$
	We can find that
	\begin{align*}
		F_2(l,a) = &~192l(l-2)(2l-1)a^4-192l(l^3 - 8l^2 + 13l - 5)a^3\\
		&-16l(l - 2)(l - 1)( l^3 + 6l^2 + 29l - 27)a^2\\
		&-16l( l^5 + 3l^4 + 4l^3 - 42l^2 + 52l - 21)a\\
		&-4l(l^5 + 3l^4 + l^3 - 24l^2 + 28l - 12).
	\end{align*}
	
	Since the coefficient of $a^4$ in $F_2(l,a)$ is positive, the expression $F_2^{(1)}(l,a)$ obtained from $F_2(l,a)$ by replacing $a^4$ with $Ka^3$ is less than $F_2(l,a)$. Then the coefficient of $a^3$ in $F_2^{(1)}(l,a)$ is $$6l(8l^2(l-3)(4l-5)+210l^2 - 395l + 142)>0.$$
	So, the expression $F_2^{(2)}(l,a)$ obtained from $F_2^{(1)}(l,a)$ by replacing $a^3$ with $Ka^2$ is less than $F_2^{(1)}(l,a)$. Then the coefficient of $a^2$ in $F_2^{(1)}(l,a)$ is $$16l^5(l-3)(6l-13)+2l^3(l-3)(285l - 98)+\frac{l}{16}(28830l^2 - 30031l + 9990)>0.$$
	Similarly, define $F_2^{(3)}(l,a)$ as the expression obtained from $F_2^{(2)}(l,a)$ by replacing $a^2$ with $Ka$. Then the coefficient of $a$ in $F_2^{(3)}(l,a)$ is
	\begin{align*}
		C = &~4l^7(l - 3)(12l - 35)+\frac{3}{4}l^5(l-3)(448l - 359)+\frac{3}{64}l^4(15320l - 29069)\\
		&+\frac{1}{512}l^2(520134l - 275585)+\frac{41061}{256}l.
	\end{align*}
	It can be seen that $C>0$ for $l\geq 3$. Finally, we have
	\begin{align*}
		&F_2^{(4)}(l,a)\\
		:=&~CK -4l(l^5 + 3l^4 + l^3 - 24l^2 + 28l - 12)\\
		=&~8l^9(l-3)(3l-11)+\frac{3}{4}l^7(l-3)(276l - 295)+\frac{1}{32}l^6( 11697l - 28745)\\
		&+\frac{3}{512}l^4(98606l - 46225)+\frac{49}{16384}l^2(30870l - 6943)+\frac{23667}{8192}l>0.
	\end{align*}
	Therefore, $$F_2(l,a)>F_2^{(1)}(l,a)>F_2^{(2)}(l,a)>F_2^{(3)}(l,a)>F_2^{(4)}(l,a)>0.\qedhere$$
\end{proof}

Here is our main result in this section.

\begin{theorem}\label{thm:spider}
	For integers $a\geq 2$ and $b\geq 2$, $$\frac{b}{2}B(\mathcal{P}_{2a+1})\leq B(\mathcal{S}_{a,b}) \leq \left(\frac{1}{2}\ln(3a+1)+T_1^\infty+\gamma-2\right)ab+o(a)b.$$
	where $\gamma$ is Euler's constant and $T_1^{\infty}$ is the absolutely convergent sum
	\begin{equation*}
		T_1^{\infty}=    
		\sum_{l=3}^{\infty}  \left(\frac{l}{l^2-l+1}-\frac{1}{l} \right)  .
	\end{equation*}
\end{theorem}

\begin{remark}
	The value of $T_1^{\infty}$ is approximately $0.3701$. Thus $T_1^\infty+\gamma-2$ is approximately $-1.0527$.
\end{remark}

\noindent
\textbf{Proof of Theorem~\ref{thm:spider}:}
From Lemma~\ref{lem:Braess on spider}, we see that each branch of $\mathcal{S}_{a,b}$ at $v$, where $v$ is the center vertex of $\mathcal{S}_{a,b}$, has $\frac{1}{2}B(\mathcal{P}_{2a+1})$ Braess edges. So, we attain a lower bound on $B(\mathcal{S}_{a,b})$.

Now, we shall find an upper bound on the number of Braess edges on $\mathcal{S}_{a,b}$. Suppose that $s\geq 0$ and $s+l\leq a$. Let $i = s+1$ and $j = s+l$. Note that $v = a+1$. Since $j<v$, the moment $\mu(\cP_{2a+1}^{\{ i,j \}},v)$ is given by \eqref{eq:moment2}. We can find that
\begin{align*}
	H_2(s,l,a):=&~\mu(\cP_{2a+1}^{\{ i,j \}},v)- \kappa\left(\cP_{2a+1}\right)-\frac{2}{3}a^2+\frac{1}{6}\\
	=&~\frac{-6(l^2 - l +1)s + 6al - 2l^3+ 2l}{3l}.
\end{align*}
We denote by $\Gamma(a)$ the number of pairs $(s,l)$ such that $H_2(s,l,a)>0$. Since $$H_2(s,l,a) =\left(\kappa(\cP_{2a+1}^{\{ i,j \}}) - \kappa\left(\cP_{2a+1}\right)\right)+\left(\mu(\cP_{2a+1}^{\{ i,j \}},v)-\kappa(\cP_{2a+1}^{\{ i,j \}})-\frac{2}{3}a^2+\frac{1}{6}\right),$$ 
it follows from Lemma~\ref{lem:Braess on spider} that $\frac{1}{2}B(\mathcal{P}_{2a+1})\leq \Gamma(a)$. We can see from \eqref{exp:spider} that for any non-edge $\{ i,j \}$, $\kappa(\mathcal{S}_{a,b}^{\{ i,j \}})-\kappa(\mathcal{S}_{a,b})$ can be evaluated by Kemeny's constant and the moment of the path $\mathcal{P}_{2a+1}$ in terms of $a$, $b$, $s$, and $l$. Moreover, recall \eqref{tempeqn:11} that since $0\leq \frac{a(b-2)}{ab+1}<1$, we have $\kappa(\mathcal{S}_{a,b}^{\{ i,j \}})-\kappa(\mathcal{S}_{a,b})<H_2(s,l,a)$. Hence, if $\{ i,j \}$ is a Braess edge, then $H_2(s,l,a)>0$. Consequently, the number of Braess edges on each branch of $S_{a,b}$ at $v$ is less than or equal to $\Gamma(a)$ and thus
\begin{align*}
	B(\mathcal{S}_{a,b}) \leq b\Gamma(a).
\end{align*}

The inequality $H_2(s,l,a)>0$ is equivalent to \begin{align*}
	0\leq s<\frac{6al-2l^3+2l}{6(l^2-l+1)}=\frac{(3a+1)l}{3(l^2-l+1)}+\frac{1}{3(l^2-l+1)}-\frac{l+1}{3}=:X(l). 
\end{align*} 
Then $6al-2l^3+2l$ is positive for $s$ to be feasible. This inequality yields 
$$l< \sqrt{3a+1} =:U_0.$$ Hence,
\begin{align*}
	\Gamma(a) = \sum_{l = 3}^{\lceil U_0\rceil -1} \lceil X(l)\rceil .
\end{align*}

Now we consider the asymptotic behaviour of $\Gamma(a)$ as $a\rightarrow \infty$. Then we see that


\begin{equation*}
	\frac{\Gamma(a)}{a}-\ln\left( U_0\right) = T_0(a)+T_1(a)+T_2(a)+T_3(a)
\end{equation*}
where
\begin{align*}
	\nonumber
	T_0(a) = & ~
	\frac{1}{a}\sum_{l=3}^{\lceil U_0\rceil -1} \left( \lceil X(l)\rceil - X(l) \right),\\
	T_1(a)  =& ~
	\sum_{l=3}^{\lceil U_0\rceil -1}\left(\frac{(3a+1)l}{3a(l^2-l+1)}-\frac{1}{l}\right) = \sum_{l=3}^{\lceil U_0\rceil -1}\left(\frac{l}{3a(l^2-l+1)}+\frac{l-1}{l(l^2-l+1)}\right),
	\\
	T_2(a) = & ~
	\frac{1}{a}\sum_{l=3}^{\lceil U_0\rceil -1} \left(\frac{1}{3(l^2-l+1)}-\frac{l+1}{3}\right),\\
	T_3(a) = & ~ -\ln\left( U_0\right)+\ln(\lceil U_0\rceil -1)+\left(-\ln(\lceil U_0\rceil -1)+\sum_{l=1}^{\lceil U_0\rceil -1}\frac{1}{l}\right)-\sum_{l=1}^{2}\frac{1}{l}.
\end{align*}  
As done in the previous section for finding the asymptotic behaviour of the number of Braess edges on path, one can find $$\lim_{a\rightarrow\infty}T_0(a)=0,\;\;\;\lim_{a\rightarrow\infty}T_3(a)=\gamma-\frac{3}{2}.$$
We can see that $$\lim_{a\rightarrow\infty}T_2(a) = -\frac{1}{2}.$$
Since $\frac{l-1}{l(l^2-l+1)}$ is a term of a convergent series, $T_1^{\infty}$ is absolutely convergent. It follows that $$\lim_{a\rightarrow\infty}T_1(a) = T_1^\infty.$$
This proves the theorem.
\hfill $\Box$

A natural question arises: Given a function $f$ such that $f(n)\leq n^2/2$ and $f(n)\rightarrow\infty$ as $n\rightarrow\infty$, can we identify a family of trees $G$ such that $B(G)\sim f$? Considering variations on spider graphs in which the lengths of the paths coming off the central vertex can vary, we are able to determine the following special case to answer the question.

\begin{theorem}\label{thm:spider2}
	For $b_1, b_2\geq 1$, $B(S_{1, b_1}\bigoplus_{v_1, v_2} S_{2, b_2}) = \binom{b_1}{2}$, where $v_1$ is the vertex of degree $b_1$ and $v_2$ is the vertex of degree $b_2$. 
\end{theorem}
\begin{proof}
	By \cite[Theorem~2.2]{ciardo2020braess}, any edge added in $S_{1, b_1}$ is Braess. Thus $B(S_{1, b_1}\bigoplus_{v_1, v_2} S_{2, b_2}) \geq \binom{b_1}{2}$. 
	Through repeated use of Theorem \ref{thm:Braess with cut vertex} it can be verified through straightforward computation that no other edges will be Braess and we have the result.
	There are 6 types of edges $\{i,j\}$ to consider.
	\begin{enumerate}[label=(\roman*)]
		\item $i,j$ are pendent vertices with $i,j\in S_{2,b_2}$
		\item $i,j\in S_{2,b_2}$ with $i$ a pendent vertex and $j$ of degree 2
		\item $i,j\in S_{2,b_2}$ with $i,j$ both of degree 2
		\item $i,j\in S_{2,b_2}$ with $i$ a pendent vertex and $j = v_2$
		\item $i,j$ are pendent vertices with $i\in S_{1, b_1}$ and $j\in S_{2,b_2}$
		\item $i$ a pendent vertex with $i\in S_{1,b_1}$ and $j\in S_{2,b_2}$ with $j$ of degree 2.
	\end{enumerate}
	Using Proposition \ref{Prop:dFd with a cut-vertex}, one can find that for $a_1\geq 1$ and $a_2\geq 1$,
	\begin{align*}
		\kappa\left(S_{1,a_1}\bigoplus\nolimits_{v_1, v_2}S_{2, a_2}\right) =&~ \frac{a_1^2 + 8a_2^2 - \frac{1}{2}a_1 - 5a_2 + 6a_1a_2}{a_1 + 2a_2},\\ 
		\mu\left(S_{1, a_1}\bigoplus\nolimits_{v_1, v_2}S_{2, a_2}, v\right) =&~ a_1 + 4a_2,
	\end{align*}
 where $v$ is the vertex of degree $a_1+a_2$.
	For the remainder of this proof we follow the notation of Theorem \ref{thm:Braess with cut vertex}.  The graphs $G_1$, $\tilde G_1$ and $G_2$ will vary depending on each case. 
	For cases (i) - (iii) $G_1 = \mathcal{P}_5$ and $G_2 = S_{1, b_1}\bigoplus_{v_1, v_2}S_{2, b_2-2}$ so by Example \ref{moment:P_n}
	\begin{align*}
		\kappa(G_1) = \frac{11}{2}, \;\;\; \mu(G_1, v) = 8.
	\end{align*}
	Case (iv) has $G_1 = \mathcal{P}_3$ and $G_2=S_{1,b_1}\bigoplus\nolimits_{v_1, v_2}S_{2, b_2-1}$, hence
	\begin{align*}
		\kappa(G_1) = \frac{3}{2},\;\;\; \mu(G_1, v) = 4.
	\end{align*}
	For cases (v) and (vi) $G_1 = \mathcal{P}_4$ and $G_2=S_{1,b_1-1}\bigoplus\nolimits_{v_1, v_2}S_{2, b_2-1}$, hence
	\begin{align*}
		\kappa(G_1) = \frac{19}{6}, \;\;\; \mu(G_1, v) = 5.
	\end{align*}
	We now work on the individual cases. Values of $\kappa(\tilde G_1)$ and $\mu(\tilde G_1, v)$ can be obtained through \eqref{eq:k(P^ij)} and Lemma \ref{lem:f_pq}
 , respectively. In cases (i) - (iii) assume that $b_2 \geq 2$, otherwise this edge does not exist. In these cases we have
 {\small
  \begin{eqnarray*}
     \mathrm{(i)}& \kappa(\tilde G_1) = 4, \mu(\tilde G_1, v) = 8,& \kappa(\tilde G) - \kappa(G) = -\frac{19b_1 + 30b_2}{2(b_1 + 2b_2)(b_1 + 2b_2 + 1)} < 0,\\
     \mathrm{(ii)}&\kappa(\tilde G_1) = \frac{39}{10}, \mu(\tilde G_1, v) = \frac{15}{2},& \kappa(\tilde G) - \kappa(G) = -\frac{b_1^2 + 4 b_1 (4 + b_2) + 4 b_2 (6 + b_2)}{
			2 (b_1 + 2 b_2) (b_1 + 2 b_2 + 1))}<0,\\
     \mathrm{(iii)}&\kappa(\tilde G_1) = \frac{59}{15}, \mu(\tilde G_1, v) = \frac{22}{3},& \kappa(\tilde G) - \kappa(G) =  -\frac{4 b_1^2 + b_1 (43 + 16 b_2) + 2 b_2 (31 + 8 b_2)}{
			6 (b_1 + 2 b_2) (b_1 + 2 b_2 + 1)} < 0.
 \end{eqnarray*}
 }
	For the remaining cases we may again assume that $b_2\geq 1$. We have
 {\small
 \begin{eqnarray*}
     \mathrm{(iv)}& \kappa(\tilde G_1) = \frac{4}{3},\;\;\;\mu(\tilde G_1, v) = \frac{8}{3}, & \kappa(\tilde G) - \kappa(G) = -\frac{-8b_1^2 - 32b_2^2 - 32b_1b_2 + b_1 + 26b_2}{6 (b_1 + 2 b_2) (b_1 + 2 b_2+1)} < 0,\\
     \mathrm{(v)}&\kappa(\tilde G_1) = \frac{5}{2},\;\;\; \mu(\tilde G_1, v) = 5, & \kappa(\tilde G) - \kappa(G) = -\frac{4 (b_1 + b_2)}{(b_1 + 2 b_2) (b_1 + 2 b_2 + 1)} < 0,\\
     \mathrm{(vi)}&\kappa(\tilde G_1) =  \frac{61}{24},\;\;\; \mu(\tilde G_1, v) = 5, & \kappa(\tilde G) - \kappa(G) = -\frac{23 b_1 + 22 b_2}{6 (b_1 + 2 b_2) (b_1 + 2 b_2+1)} < 0.
 \end{eqnarray*}
 }\qedhere
\end{proof}

\section{Braess edges on a broom}\label{sec:broom}
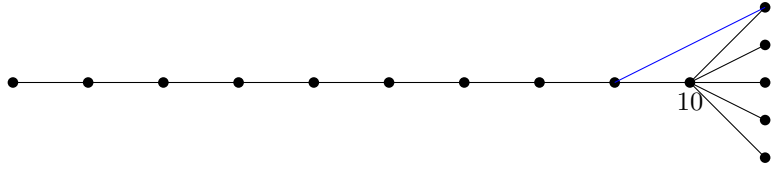
\begin{figure}[h]
	\centering
	\begin{tikzpicture}
		\tikzstyle{every node}=[circle,  fill=black, minimum width=4pt, inner sep=1pt]
		\draw{
			(0,0)node{}--(1,0)node{}--(2,0)node{}--(3,0)node{}--(4,0)node{}--(5,0)node{}--(6,0)node{}--(7,0)node{}--(8,0)node{}--(9,0)node[label={[shift={(0,-.6)}]{$10$}}]{}--(10, 1)node{}
			(9,0)--(10,0.5)node{}
			(9,0)--(10,0)node{}
			(9,0)--(10, -0.5)node{}
			(9,0)--(10,-1)node{}
		};
		\draw[blue]{(8,0)--(10,1)};
	\end{tikzpicture}
	\caption{The broom $\mathcal{B}_{10,5}$. The blue edge here creates a cycle of length $l=3$, giving $\mathcal{B}_{10, 5}^3$.}
	\label{fig:Broom}
\end{figure}

Let $\mathcal{P}_k = (1, \hdots, k)$. A \emph{broom} is the graph $\mathcal{B}_{k, p} = \mathcal{P}_k \bigoplus_{k,v} S_{1, p}$ with $p\geq 2$, where $v$ is the vertex in $\mathcal{S}_{1,p}$ with degree $p$. See Figure \ref{fig:Broom} for an example. We let $B_1$ be the number of Braess edges joining two vertices of $S_{1, p}$ for $\mathcal{B}_{k, p}$, $B_2$ be the number of Braess edges for $\mathcal{B}_{k, p}$ joining a vertex of $\mathcal{P}_k$ and a pendent vertex of $S_{1, p}$, and $B_3$ be the number of Braess edges on $\mathcal{P}_k$ for $\mathcal{B}_{k, p}$. Then, $B(\mathcal{B}_{k,p}) = B_1+B_2+B_3$. Since any non-edge between twin pendent vertices is a Braess edge by \cite[Theorem~2.2]{ciardo2020braess}, we have $B_1=\frac{p(p-1)}{2}$. We shall investigate $B_2$ and $B_3$ in Subsections 5.1 and 5.2, respectively, with the information about where Braess edges occur.

Here is the main result of this section.
\begin{theorem}\label{thm:broom}
	Let $B_1$, $B_2$ and $B_3$ be defined as above. Then,
	\begin{align*}
		B_1 = \frac{p(p-1)}{2}, \;\;\;B_2 = O(k),\;\;\;B_3 = \Theta(k\ln(k)).
	\end{align*}
	This implies that as $k\rightarrow \infty$, we have 
	\[\begin{array}{cll}
		B(\mathcal{B}_{k,p})\;\; =& \Theta(k \ln (k)+p^2), & \text{if $p/k\rightarrow0$;}\\
		B(\mathcal{B}_{k,p})\;\; \sim& \frac{1}{2}p^2, & \text{if $\frac{p}{k}\rightarrow\beta\in (0,\infty].$}
	\end{array}
	\]
\end{theorem}
\begin{proof}
	As we shall see below, in all cases that $p/k$ tends to a finite number or infinity, Corollary~\ref{cor:B2isO(k)} tells us that $B_2 = O(k)$, and Theorem~\ref{thm:B3} provides $B_3 = \Theta(k \ln(k))$.
\end{proof}

\begin{remark}
	We see a threshold at $p=\sqrt{k \ln(k)}$:  If $p$ grows faster than this rate, then $p^2$ dominates and most of the Braess edges for the broom occur on $\mathcal{S}_{1,p}$; if $p$ grows more slowly than this rate, then $k\ln(k)$ dominates, and most of the Braess edges occur on $\mathcal{P}_k$.
\end{remark}

\subsection{Braess edges from pendent to path}

Proposition \ref{prop:KemDifferenceLollipopBroom} describes the change in Kemeny's constant in a broom $\mathcal{B}_{k,p}$ after adding an edge from one of the $p$ pendent vertices in $S_{1, p}$ to one of the vertices in $\mathcal{P}_k$ to create a cycle of length $l$. Denote the resulting graph by $\mathcal{B}_{k,p}^l$. An edge of this sort will go across the 1-separation suggested in the definition of $\mathcal{B}_{k,p}$. In order to still make use of the 1-separation formulas for Kemeny's constant we notice we can write $\mathcal{B}_{k, p} = \mathcal{P}_{k+1}\bigoplus_{k,v} S_{1, p-1}$. Thinking of the graph in this way, the edge added to get $\mathcal{B}_{k,p}^l$ will not go across the 1-separation, but be entirely contained in the $\mathcal{P}_{k+1}$ part.

In considering Braess edges of this sort in a Broom $\mathcal{B}_{k,p}$, it will be helpful to understand Kemeny's constant for a lollipop graph, $\mathcal{L}_{l, n}$. Here we recall that this is a graph obtained by taking a path $\mathcal{P}_n = (1, \hdots, n)$ and adding an edge $(n, n-l+1)$ to create a cycle of length $l$ at one end of the graph. Notice that $\mathcal{B}_{k,p}^l = \mathcal{L}_{l, k+1}\bigoplus_{k, v} S_{1, p-1}$.
\begin{lemma}\label{lem:LollipopValues}
	Let $k\geq 3$ and $3\leq l\leq k+1$. Then    \begin{align*}
		\kappa(\mathcal{L}_{l, k+1}) &= \frac{2(k+1)^3 - 4(k+1)l^2 + 3l^3 - (k+1)}{6(k+1)}\\
		\mu(\mathcal{L}_{l, k+1}, k) &= \frac{(l+1)(l-1)}{3} + (k - l + 1)^2 + \frac{4(k - l + 1)(l-2)}{l}.
	\end{align*}
\end{lemma}
\begin{proof}
	The formula for Kemeny's constant was given in (3.4) of \cite{kirkland2016kemeny}. The result on the moment comes from using \cite[Lemma~2.3]{faught20221}.
\end{proof}

Though Theorem \ref{thm:Braess with cut vertex} already gives an expression for the change in Kemeny's constant, here we give another. The following can be found by using Theorem 2.1 of \cite{faught20221} together with
Examples \ref{moment:P_n} and \ref{moment:S_n} and Lemma \ref{lem:LollipopValues} to find $\kappa(\tilde G)$ and $\kappa(G)$, then taking the difference.

\begin{proposition}\label{prop:KemDifferenceLollipopBroom}
	Let $k\geq 3$, $p\geq 2$ and $3\leq l\leq k+1$. Then 
	\begin{align*}
		&\kappa(\mathcal{B}_{k, p}^l) - \kappa(\mathcal{B}_{k,p})\\
		=& \left[\frac{k+1}{k+p}\kappa(\mathcal{L}_{l,k+1}) - \frac{k}{k+p-1}\kappa(\mathcal{P}_{k+1})\right] +\left[ \frac{p-1}{k+p}\mu(\mathcal{L}_{l,k+1}, k) - \frac{p-1}{k+p-1}(k^2-2k+2)\right]\\
		&\hspace{4em}+ \frac{p-1}{2(k+p)(k+p-1)}.
	\end{align*}
\end{proposition}

These next three lemmas will be helpful in ruling out edges that will not be Braess edges.

\begin{lemma}\label{lem:KemDifferenceDecreasingInl}
	Let $k\geq 5$ and $p\geq 2$. Then $\kappa(\mathcal{B}_{k, p}^l) - \kappa(\mathcal{B}_{k,p})$ is decreasing in $l$ for $3\leq l \leq (k+1)/2$.
\end{lemma}
\begin{proof}
	Note that only two terms 
	in Proposition \ref{prop:KemDifferenceLollipopBroom} depend on $l$. $\kappa(\mathcal{L}_{l, k+1})$  was observed to be decreasing in $l$ for $3\leq l\leq 8(k+1)/9$ in Example 3.4 of \cite{kirkland2016kemeny}.
	
	Now we show that $\mu(\mathcal{L}_{l, k+1}, k)$ is decreasing in $l$ on $3\leq l \leq (k+1)/2.$
	Notice
	\[
	\frac{d}{dl}\left[\mu(\mathcal{L}_{l, k+1},k)\right] = \frac{8}{3}l - 2k - 6 + \frac{8(k+1)}{l^2}.
	\]
	Taking the second and third derivative we see 
	\begin{align*}
		\frac{d^2}{dl^2}\left[\mu(\mathcal{L}_{l, k+1},k)\right] =&\: \frac{8}{3} - \frac{16(k+1)}{l^3}\\
		\frac{d^3}{dl^3}\left[\mu(\mathcal{L}_{l, k+1},k)\right] = &\:\frac{48(k+1)}{l^4}.
	\end{align*}
	Note that the third derivative is clearly positive for all positive values of $k$,
	hence the first derivative is a convex function. Thus if we can find two values of $l$ which give a negative value for the first derivative, then all intermediate values of $l$ will also give a negative value for the first derivative. Plugging in $l=3$ and $l= (k+1)/2$ gives
	\begin{align*}
		\left.\frac{d}{dl}\left[\mu(\mathcal{L}_{l, k+1},k)\right]\right\rvert_{l=3} &= \frac{26}{9} - \frac{10k}{9} < 0 &\text{if $k\geq 3$.}\\
		\left.\frac{d}{dl}\left[\mu(\mathcal{L}_{l, k+1},k)\right]\right\rvert_{l=(k+1)/2} &= \frac{-2(k^2 + 8k - 41)}{3(k+1)} < 0 &\text{if $k\geq 4$.}
	\end{align*}
	Thus we have the desired conclusion.
\end{proof}

\begin{lemma}\label{lem:2ndTermNegative}
	Let $p\geq 2$ and $k\geq 5$. Then 
	\[
	\mu(\mathcal{L}_{l, k+1}, k) - (k^2 - 2k + 2) < 0.
	\]
\end{lemma}
\begin{proof}    
	We use some of the calculations from the proof of Lemma \ref{lem:KemDifferenceDecreasingInl}. For a given $k\ge5$, let $p(l) = \mu(\mathcal{L}_{l, k+1}, k) - (k^2 - 2k + 2)$.  Since $\left.\frac{d}{dl}p(l)\right\rvert_{l=3}<0$ and $\frac{d}{dl}p(l)$ is a convex function, it suffices to check both $p(3)<0$ and $p(k+1)<0$. It is straightforward to verify that
	\begin{align*}
		\mu(\mathcal{L}_{3, k+1}, k) - (k^2 - 2k + 2) &= 2 - \frac{2k}{3} < 0, \\
		\mu(\mathcal{L}_{k+1, k+1}, k) - (k^2 - 2k + 2) &= -\frac{2}{3}(k-1)(k-3) < 0. \qedhere
	\end{align*}
\end{proof}
Note that Lemma \ref{lem:2ndTermNegative} implies that the second term in square brackets in Proposition \ref{prop:KemDifferenceLollipopBroom} is strictly negative.

\begin{lemma}\label{lem:notBraessPathGivesNotBraessBroom}
	Let $k\geq5$. If the edge added into $\mathcal{B}_{k,p}$ to get $\mathcal{B}_{k,p}^l$ is not a Braess edge for $\mathcal{P}_{k+1}$, then it is not a Braess edge for $\mathcal{B}_{k,p}$.
\end{lemma}
\begin{proof}
	The expression in the second square bracket in Proposition \ref{prop:KemDifferenceLollipopBroom} is always negative as was shown in Lemma \ref{lem:2ndTermNegative}. Further, we see that the maximum is attained at $l=3$ or at $l = k + 1$. We examine these two cases. Let $l = 3$. This gives 
	\[
	\left[ \frac{p-1}{k+p}\mu(\mathcal{L}_{3, k+1}, k) - \frac{p-1}{k+p-1}(k^2-2k+2)\right] = -\frac{(p-1)(5k^2 + 2k(p-7) - 6(p-2))}{3(k+p)(k+p-1)}.
	\]
	
	Now under the assumption that the edge is not a Braess edge, we have $\kappa(\mathcal{L}_{l, k+1}) \leq \kappa(\mathcal{P}_{k+1}) = (2k^2 + 1)/6$. Then
	\begin{align*}
		\frac{k+1}{k+p}\kappa(\mathcal{L}_{l, k+1}) - \frac{k}{k+p-1}\kappa(\mathcal{P}_{k+1}) &\leq \frac{2k^2+1}{6}\left(\frac{k+1}{k+p} - \frac{k}{k+p-1}\right)\\
		&= \frac{(2k^2 + 1)(p-1)}{6(k+p)(k+p-1)}.
	\end{align*}
	Simplifying the expression in Proposition \ref{prop:KemDifferenceLollipopBroom} gives 
	\begin{align*}
		&\kappa(\mathcal{B}_{k,p}^3) - \kappa(\mathcal{B}_{k,p})\\
		&\leq \frac{(2k^2+1)(p-1)}{6(k+p)(k+p-1)} - \frac{(p-1)(5k^2+2k(p-7)-6(p-2))}{3(k+p)(k+p-1)} +\frac{p-1}{2(k+p)(k+p-1)}\\
		&= -\frac{2(p-1)}{3(k+p)(k+p-1)}[2k^2 + k(p-7) - 3p + 5].
	\end{align*}
	This is negative if $p\geq 2$. 
	
	Now suppose that $l = k+1$. Similar work as above gives
	\begin{align*}
		\kappa(\mathcal{B}_{k,p}^{k+1}) - \kappa(\mathcal{B}_{k,p}) &\leq -\frac{2(k-1)(p-1)(k^2+k(p-3)-3p+1)}{3(k+p)(k+p-1)}.
	\end{align*}
	This is easily seen to be negative if $p\geq 2$. 
	
	Therefore, if a non-edge is not Braess in $\mathcal{P}_{k+1}$ it is not Braess in $G$.
\end{proof}

The following proposition is the entire formula in the Proposition \ref{prop:KemDifferenceLollipopBroom}, expanded out. We include this version too since different insights are readily seen looking at the expression in this manner than are seen from Proposition \ref{prop:KemDifferenceLollipopBroom}.

\begin{proposition}\label{prop:KemDifferenceLollipopBroomExpanded}
	Let $k\geq 3$, $p\geq 2$ and $3\leq l\leq k+1$. Then
	\begin{align*}
		&\kappa(\mathcal{B}_{k,p}^l) - \kappa(\mathcal{B}_{k,p})\\
		=&\left(\frac{1}{6l(k+p-1)(k+p)}\right)\bigg(4k^3l - 4k^2[l^3 + 3l^2(p-1) - l(12p-11)+12(p-1) ]\\
		&\quad+ k[3l^4 + 4l^3(p-2)  - 12l^2(p^2 + p-2) + 16l(3p^2-p-2) - 48(p-1)p]\\
		&\quad+ (l-2)(p-1)[3l^3 + 24(p-1) + 2l^2(4p-3)-4l(5p-6)]\bigg).
	\end{align*}
\end{proposition}
Looking at Proposition \ref{prop:KemDifferenceLollipopBroomExpanded}, it is clear that for fixed $l$ and $p$, if the path is long enough ($k$ is large enough) then the cycle is Braess. While we do not pin down exactly how large $k$ must be to ensure an edge is Braess, these next few results partially answer that question. In fact, Proposition \ref{prop:KHowBigBraessPendantToPath} will give a sufficient size of $k$ for the edge to be Braess.

\begin{lemma}\label{lem:sqrtK+1NotBraess}
	Let $k\geq 6$ and $p\geq 2$. If $l = \sqrt{k} + 1$, then $\kappa(\mathcal{B}_{k,p}^l) - \kappa(\mathcal{B}_{k,p}) < 0 $.
\end{lemma}
\begin{proof}
	As the denominator in Proposition \ref{prop:KemDifferenceLollipopBroomExpanded} is positive, we see that $\kappa(\mathcal{B}_{k,p}^l) - \kappa(\mathcal{B}_{k,p}) < 0 $ if and only if the numerator is negative. Then putting $l = \sqrt{k} + 1$ in the numerator gives that $\kappa(\mathcal{B}_{k,p}^l) - \kappa(\mathcal{B}_{k,p}) < 0$ if and only if
	\begin{align*}
		-(\sqrt{k}+1)[k^{5/2}(12p-7) - k^2(16p-21) + k^{3/2}(12p^2 - 7p - 6)&\\
		- k(20p^2 - 29p + 10) + \sqrt{k}(4p^2 - 25p + 21) - 3(4p^2 - 5p + 1)] &< 0.
	\end{align*}
	
	This is true if and only if the expression in the square bracket is positive.
	Analyzing the $k^{5/2}$ and $k^2$ terms together, they can be seen to be positive for $k\geq 4.$ The $\sqrt{k}$ term is positive if $p\geq 6$. Looking at the $k^{3/2}, k, $ and $k^0$ terms together these are also seen to be positive if $k\geq 4$. This is seen by replacing the constant 3 in the $k^0$ term with $k$. 
	
	Checking the $2\leq p\leq 5$ cases one finds that the expression in the square bracket is positive if $k\geq 6$. The conclusion follows.
\end{proof}

\begin{proposition}\label{prop:NotBraessInBroom}
	Let $p\geq 2$ and $k\geq 6$. Suppose we add an edge to $\mathcal{B}_{k,p}$ to obtain $\mathcal{B}_{k,p}^l$. If $l \geq \sqrt{k} + 1$ then the added edge is not Braess. This implies that $B_2 = O(p\sqrt{k})$ as $k\rightarrow \infty$.
\end{proposition}
\begin{proof}
	By Lemma \ref{lem:sqrtK+1NotBraess} $l = \sqrt{k} + 1$ does not give a Braess edge. By Lemma \ref{lem:KemDifferenceDecreasingInl} there can be no Braess edges of this sort for $l\in[\sqrt{k}+1, (k+1)/2]$. Further, by Theorem \ref{thm:nonedge crossing the center} combined with Lemma \ref{lem:notBraessPathGivesNotBraessBroom} there are no Braess edges for $l\geq (k+1)/2$.
\end{proof}

\begin{corollary}\label{Cor:l implies s for no Braess}
	Let $p\geq 2$ and $k\geq 6$. If  $\kappa(\mathcal{B}_{k,p}^l) - \kappa(\mathcal{B}_{k,p})<0$ for some $l\geq 3$, then $\kappa(\mathcal{B}_{k,p}^s) - \kappa(\mathcal{B}_{k,p})<0$ for all $l\leq s\leq k+1$.
\end{corollary}
\begin{proof}
	We note that $\sqrt{k}+1<\frac{k+1}{2}$ for $k\geq 6$. If $l\geq \sqrt{k}+1$ then the result follows from Proposition~\ref{prop:NotBraessInBroom}. If $l< \sqrt{k}+1$  then the conclusion follows from Lemma~\ref{lem:KemDifferenceDecreasingInl} and Proposition~\ref{prop:NotBraessInBroom}.
\end{proof}

\begin{proposition}\label{prop:KHowBigBraessPendantToPath}
	Let $l\geq 3$ and $p\geq 2$. Suppose we add an edge to $\mathcal{B}_{k,p}$ to obtain $\mathcal{B}_{k,p}^l$. If $k \geq (0.007p + 1)l^2+3(p-1)l - 7(p-1)$, then the edge added is Braess.
\end{proposition}
\begin{proof}
	As seen in Proposition \ref{prop:KemDifferenceLollipopBroomExpanded}, the numerator of $\kappa(\mathcal{B}_{k,p}^l) - \kappa(\mathcal{B}_{k,p})$ is a cubic in $k$, hence there will come a point for which $\kappa(\mathcal{B}_{k,p}^l) - \kappa(\mathcal{B}_{k,p})$ begins increasing in $k$ indefinitely. We first find that point, then show that the estimate for $k$ given in this Proposition statement is past that, then show that $\kappa(\mathcal{B}_{k,p}^l) - \kappa(\mathcal{B}_{k,p}) > 0$ at that estimate. This will show that if $k\geq K_1 = (0.007p + 1)l^2+3(p-1)l - 7(p-1)$ then the edge added will be Braess.
	
	Using Proposition \ref{prop:KemDifferenceLollipopBroomExpanded}, we take the derivative with respect to $k$ of the numerator. This has roots
	\begin{equation}\label{eq:kIncreasingPastHere}
		k = g_1(l,p) \pm \frac{\sqrt{g_2(l, p)}}{24l},
	\end{equation}
	where
	\begin{align*}
		g_1(l, p) &~= \frac{1}{3}l^2 + l(p-1) - \frac{1}{3}(12p-11) + \frac{4(p-1)}{l},\\
		g_2(l, p) &~= 64(l^3 + 3l^2(p-1) - l(12p-11) + 12(p-1))^2\\
		&\quad - 48l(3l^4 + 4l^3(p-2) - 12l^2(p^2 + p - 2) + 16l(3p^2 - p - 2) - 48(p-1)p).
	\end{align*}
	Let $K_2$ be the right hand side of equation (\ref{eq:kIncreasingPastHere}) with the positive square root. It can be shown that $g_2(l,p)\geq 0$ for all $p\geq 2$ and $l\geq 3$, hence $K_2$ is a real number. We want to show that $K_1 - K_2 \geq 0$. Straightforward algebra gives this is true if and only if
	\begin{align*}
		&(24l)^2(K_1-g_1(l,p))^2-g_2(l,p)\\
		=&~\frac{3}{15625}( 7p + 1000)(21p + 1000)l^6 + \frac{144}{125}(14p^2 + 986p - 875)l^5 \\
		&+ \frac{48}{125}( 4437p^2 - 10430p + 6500)l^4 - \frac{2304}{125}(158p + 125)(p - 1)l^3 \\
		&-576(p - 1)(27p - 29)l^2 + 2304(p - 1)(13p - 14)l\geq 0.
	\end{align*}
	
	Looking at the $l^6, l^4, $ and $l^2$ terms together it can be seen that these three terms are positive if $l\geq 4$ and $p\geq 2$. Looking at $l^5$ and $l^3$  terms together they are seen to be positive if $l \geq 14$ and $p\geq 2$.
	The $l$ term is positive provided $p\geq 2$. Checking when $l=3, 4, \hdots, 13$ it is seen that the entire expression is positive if $p\geq 2$ and $l\geq 3$. Therefore $K_1 - K_2\geq 0$ for $p\geq 2$ and $l\geq 3$.
	
	Now as $\kappa(\mathcal{B}_{k,p}^l) - \kappa(\mathcal{B}_{k,p})$ is increasing in $k$ for $k\geq K_1$, we show that this difference is indeed positive at $k = K_1$. Plugging $K_1$ in for $k$ and then doing work similar to what was done just above shows that $\kappa(\mathcal{B}_{k,p}^l) - \kappa(\mathcal{B}_{k,p})\geq 0$ at $k = K_1$ and we have the result.
\end{proof}

Recall that $B_2$ denotes the number of Braess edges for $\mathcal{B}_{k,p}$ connecting a pendent vertex in $S_{1,p}$ to a vertex in $\mathcal{P}_k$. Now we present several results regarding $B_2$.

\begin{theorem}\label{thm:pTooBigForBraess}
	If $p\geq k\geq 3$ then $B_2 =0$.
\end{theorem}
\begin{proof}
	
	We shall first show that $\kappa(\mathcal{B}_{k,p}^l) - \kappa(\mathcal{B}_{k,p})<0$ for $p\geq k\geq 6$. From Corollary~\ref{Cor:l implies s for no Braess}, it suffices to show that $\kappa(\mathcal{B}_{k,p}^3) - \kappa(\mathcal{B}_{k,p}) < 0$. Now consider
	$$\kappa(\mathcal{B}_{k,p}^3) - \kappa(\mathcal{B}_{k,p}) = \frac{-4(k-3 )p^2 - (4k^2-13)p + 4k^3 - 28k^2 + 49k - 25}{6(k + p)(k + p - 1)}.$$
	Given $k\geq 6$, the maximum of the numerator is attained at $p=k$, which is given by $-4k^3 - 16k^2 + 62k - 25<0$. Hence, $B_2 = 0$ for $p\geq k\geq 6$.
	
	It remains to show that $\kappa(\mathcal{B}_{k,p}^l) - \kappa(\mathcal{B}_{k,p}) < 0$ for $3\leq k\leq 5$ and $p\geq k$. Using the expression in Proposition~\ref{prop:KemDifferenceLollipopBroomExpanded}, computing for each of pairs $(l,k)$ with  $3\leq k\leq 5$ and $3\leq l\leq k+1$, one can establish the desired result.
\end{proof}

\begin{theorem}\label{thm:pTooBigForBraessAsymptotically}
	Fix a real number $\beta>0$.  Assume 
	$\frac{p}{k}\rightarrow \beta$ as $k\rightarrow \infty$. 
	If $\beta \geq \frac{-1 + \sqrt{5}}{2}$, then $B_2 = 0$ for all
	sufficiently large $k$. If $0 < \beta < \frac{-1 + \sqrt{5}}{2}$, 
	then $B_2 = \Theta(p) = \Theta(k)$ as $k\to\infty$.
\end{theorem}
\begin{proof}
	We assume that $\frac{p}{k}\to\beta$.
	Showing that the edge to $\mathcal{B}_{k,p}$ to obtain $\mathcal{B}_{k,p}^l$ is Braess is equivalent to showing that the numerator in Proposition \ref{prop:KemDifferenceLollipopBroomExpanded} is positive. 
	
	First, we prove the following claim:  There is a constant $C$ (depending on $\beta$) such that, for all sufficiently large $k$, no edge with $l>C$ is Braess.  Assume that the claim is false. Then there is a sequence $l=l(k)$ such that $l(k)\to\infty$ and the edge to $\mathcal{B}_{k,p}$ to obtain $\mathcal{B}_{k,p}^{l(k)}$ is Braess for $\mathcal{B}_{k,p}$ (for infinitely many values of $k$). 
	By Proposition \ref{prop:NotBraessInBroom}, we have $l<\sqrt{k}+1$.
	We then see that the dominant terms in the numerator of Proposition \ref{prop:KemDifferenceLollipopBroomExpanded} are of order 
	$k^3l^2$, specifically $-12k^2l^2p-12kl^2p^2$. Since these terms 
	are negative, and all the other terms are $o(k^3l^2)$, we have a contradiction.  This proves the claim, and shows that $B_2=O(k)$.
	
	Given the above claim, we can now assume $l=O(1)$.
	Using $p = \beta k+o(k)$ and simplifying the numerator in Proposition \ref{prop:KemDifferenceLollipopBroomExpanded}, we get
	\begin{equation}
		\label{eq.kcubednumerator}
		k^3 (-48 \beta - 48 \beta^2 + 4 l + 48 \beta l + 48 \beta^2 l - 12 \beta l^2 - 12 \beta^2 l^2) + o(k^3).
	\end{equation}
	If for a given $\beta$ and $l$ the $k^3$ coefficient is negative, then in the limit as $k\to\infty$ we know that the numerator must be negative, hence an edge with that given $l$ is not a Braess edge. Further, by the monotonicity argument of Corollary~\ref{Cor:l implies s for no Braess}, no edge with $\hat l \geq l$ will be Braess. Thus, in particular, if $l = 3$ and this coefficient is negative, then we have $B_2 = 0$ for all sufficiently large $k$.
	
	Solving for the roots of this coefficient as a quadratic in $\beta$ we see that this coefficient is negative if $\beta$ is greater than the larger of the two roots, hence if 
	\[
	\beta > -\frac{1}{2} + \frac{\sqrt{3(3l^2 - 8l + 12)}}{6(l-2)}
	\]
	If $l = 3$ then the value of the right hand side is $\frac{-1 + \sqrt{5}}{2}\approx 0.618$. Hence if $\beta > \frac{-1 + \sqrt{5}}{2}$ then $B_2 = 0$ as $k\to\infty.$
	One can check that if $\beta = \frac{-1 + \sqrt{5}}{2}$ then $B_2 = 0$ for all sufficiently large $k$ in this case too.
	Finally, if $0 < \beta < \frac{-1 + \sqrt{5}}{2}$, then Equation (\ref{eq.kcubednumerator}) is positive for $l=3$ as $k\rightarrow\infty$, showing that $B_2\geq p$, and hence that $B_2=\Theta(p)=\Theta(k)$ (since our earlier claim implies that $B_2\leq Cp$).
\end{proof}

\begin{theorem}\label{thm:BraessOrderKDominatesP}
	Let $n = p+k$ be the number of vertices in $\mathcal{B}_{k,p}$. If $\frac{p}{k}\to 0$ as $n\to\infty$, then $B_2 = O(k)$.
\end{theorem}
\begin{proof}
	From Proposition~\ref{prop:NotBraessInBroom}, we have $B_2 = O(p\sqrt{k})$ as $k\to \infty$. If $p$ is fixed, then $B_2 = O(\sqrt{k})$. We assume that $\frac{p}{k}\to0$ as $k,p\to \infty$.
	To determine if an edge is Braess, it will suffice to show that the numerator in Proposition \ref{prop:KemDifferenceLollipopBroomExpanded}, which we will denote as $F(k, p, l)$, is positive or negative. We will find an upper bound on $B_2$ by showing that as $n \to\infty$, if $l\geq k/p$, the edge is not Braess.
	
	We claim that if $p\geq 4$ then $F(k,p,l)<0$ for $k\geq 5p$ and $l\geq k/p$. Let $H(k, p) = p^4F(k, p, k/p)$. We have
	\begin{align*}
		H(k, p) =&\: k^5(3 - 4p) - k^4(p-1)(8p^2 - 8p - 3) + 4k^3p(9p^3-12p^2+p+3)\\
		&- 4k^2p^2(p-9)(p-1) - 8kp^3(p-1)(6p^2-8p+9) - 48p^4(p-1)^2.
	\end{align*}
	Let $c_i$ denote the coefficient of $k^i$. Notice that $c_5, c_4, c_1, c_0 < 0$ if $p\geq 2$. Further $c_3 > 0$ if $p\geq 2$. Regardless of the sign of $c_2$, Descartes' Sign Rule implies that there are at most 2 positive roots to $H(k, p)$. We will now use this fact to show that $H(k, p) < 0$ as $k+p\to\infty$.
	
	We can find that for $p\geq 4$,
	\begin{align*}
		H(p, p) =& -p^4(20p^3 - 24p^2 - 2p + 3)<0,\\
		H(3p, p) =& 15p^4(12p^3 - 48p^2 + 32p - 5)>0,\\
		H(5p, p) =& -p^4(740p^3 + 8088p^2 - 7166p + 963)<0.
	\end{align*}
	It follows that $H(k, p) < 0$ for $k\geq 5p$. By the monotonicity argument of Corollary~\ref{Cor:l implies s for no Braess}, if $p\geq 4$ then $F(k,p,l)<0$ for $k\geq 5p$ and $l\geq k/p$. Therefore, if $p/k\to 0$ as $k,p\to\infty$ then $B_2 < p\cdot\frac{k}{p} = k$ and hence $B_2 = O(k)$. 
\end{proof}

We summarize the previous three results in the following corollary.
\begin{corollary}\label{cor:B2isO(k)}
	Let $n = p+k$ be the number of vertices in $\mathcal{B}_{k,p}$. Then $B_2 = O(k)$ as $n\to\infty.$
\end{corollary}

\subsection{Braess edges on the path in a broom}
Let $\mathcal{P}_k = (1, \hdots, k)$. Consider the broom $\mathcal{B}_{k, p} = \mathcal{P}_k \bigoplus_{k,v} S_{1, p}$ with $k\geq 4$ and $p\geq 2$. 
Let $1\leq i<j\leq k$ with $j-i\geq 2$. We use $\mathcal{B}_{k,p}^{\{ i,j \}}$ to denote the graph obtained from $\mathcal{B}_{k,p}$ by adding the edge $\{ i,j \}$ into $\mathcal{B}_{k,p}$.

\begin{proposition} \label{prop: kemeny B-B}
	Let $k\geq 4$ and $p\geq 2$. Let $1\leq i<j\leq k$. Then,
	\begin{align*}
		&\kappa(\mathcal{B}_{k,p}^{\{ i,j \}})-\kappa(\mathcal{B}_{k,p})\\
		=&~\kappa(\cP_{k}^{\{ i,j \}}) - \kappa\left(\cP_{k}\right)+\frac{p}{k+p}\left(\mu(\mathcal{P}_k^{\{ i,j \}}, k)-\kappa(\cP_{k}^{\{ i,j \}})-\frac{1}{2}\right)-\frac{p}{k+p-1}\left(\frac{2}{3}(k-1)^2-\frac{2}{3}\right).
	\end{align*}
\end{proposition}
\begin{proof}
	From \eqref{moment:S_n}, $\mu(S_{1, p}, k)-\kappa(S_{1, p})=\frac{1}{2}$. It is straightforward from Theorem~\ref{thm:Braess with cut vertex}. 
\end{proof}

Let $l = j-i+1$ and $s = i-1$. It follows from Lemma~\ref{lem:f_pq} that
\begin{align*}
	\mathbf{d}_{\cP_{k}^{\{ i,j \}}}^TF_{\cP_{k}^{\{ i,j \}}}\be_k=lk^2 - 2(l^2 -l +1)s - \frac{1}{3}(2l^3+l).
\end{align*}
Hence, 
\begin{align} \label{eq: mu Pk}
	\mu(\mathcal{P}_k^{\{ i,j \}}, k) =k^2 - 2\left(l -1 +\frac{1}{l}\right)s - \frac{1}{3}(2l^2+1).
\end{align}
Here we recall \eqref{eq:k(P^ij)} that 
\begin{align}\label{eq: k Pij}
	\kappa(\cP_{k}^{\{ i,j \}})=& \frac{1}{6kl}( 12(l^2-l+1)s^2 + 12(l^3 - l^2k - l^2 + lk + l - k)s \\
	\notag &+ l(3l^3 -4l^2k +2k^3 -k)).
\end{align}
Let 
\[
\kappa(\mathcal{B}_{k,p}^{\{ i,j \}})-\kappa(\mathcal{B}_{k,p}) = \frac{D(s,l,k,p)}{6l(k+p-1)(k+p)}.
\]
Then $D(s,l,k,p)>0$ if and only if $\kappa(\mathcal{B}_{k,p}^{\{ i,j \}})-\kappa(\mathcal{B}_{k,p})>0$. Simplifying  $D(s,l,k,p)$ using Example~\ref{moment:P_n}, Proposition~\ref{prop: kemeny B-B}, \eqref{eq: mu Pk} and \eqref{eq: k Pij} and collecting it in terms of $s$ yields 
\begin{align*}
	D(s,l,k,p)=& ~12(k + p-1)(l^2 - l + 1)s^2 - 12(k+ p-1)(l^2 - l + 1)(k + p - l )s \\
	&~+4l(-l^2 + 3k -2)p^2 + l(12k^2 - 8kl^2 -16k + 3l^3 + 4l^2 + 8)p \\
	&~+l(k-1)(4k^2 - 4kl^2 - 4k + 3l^3).
\end{align*}
One can verify that two roots are given by
\begin{align}\label{s_11}
	s_1(l,k,p)= \frac{k+p-l}{2}-\sqrt{C(l,k,p)}\;\;\text{and}\;\; s_2(l,k,p)=\frac{k+p-l}{2}+\sqrt{C(l,k,p)},
\end{align}
where
\begin{multline}\label{eqn:c(l,k,p)}
	C(l,k,p)
	=\frac{1}{4}(k+p)^2-\frac{lkp}{l^2-l+1}-\frac{l}{2}(k+p)+\frac{l^3-l}{3(l^2-l+1)}(k+p)\\
	+\frac{lp}{l^2-l+1}-\frac{l^3-l^2}{4(l^2-l+1)}-\frac{k}{(k+p-1)}\frac{(k^2-3k+2)l}{3(l^2-l+1)}.
\end{multline}
Hence, $\{ i,j \}$ is a Braess edge for $\mathcal{B}_{k,p}$ if and only if 
\begin{align}\label{range of s}
	0\leq s<s_1(l,k,p)\;\;\text{or}\;\;s_2(l,k,p)<s\leq k-l+1.
\end{align}

We now consider where a Braess edge forms when $n$ tends to infinity, provided that $l$ is given. Setting $\beta = \frac{p}{k}$, we can find that as $k\rightarrow \infty$,
\begin{align*}
	\frac{1}{k}\left(\frac{k+p-l}{2}\pm\sqrt{C(l,k,p)}\right)\rightarrow \frac{1+\beta}{2}\pm\sqrt{\frac{(1+\beta)^2}{4}-\frac{l\beta}{l^2-l+1}-\frac{l}{3(1+\beta)(l^2-l+1)}}.
\end{align*}
Let $n= k+p$. It follows that 
$$\frac{s_1(l,k,p)}{k}\rightarrow\begin{cases*}
	\frac{l}{l^2-l+1}, & \text{if $\beta\rightarrow\infty$ as $n\rightarrow\infty$;}\\
	\frac{1}{2}-\sqrt{\frac{1}{4}-\frac{l}{3(l^2-l+1)}}, & \text{if $\beta\rightarrow 0$ as $n\rightarrow\infty$.}
\end{cases*}$$
We also have
$$\frac{s_2(l,k,p)}{k}\rightarrow\begin{cases*}
	\infty, & \text{if $\beta\rightarrow\infty$ as $n\rightarrow\infty$;}\\
	\frac{1}{2}+\sqrt{\frac{1}{4}-\frac{l}{3(l^2-l+1)}}, & \text{if $\beta\rightarrow 0$ as $n\rightarrow\infty$.}
\end{cases*}$$
Therefore, if $k$ dominates over $p$, then the location in which Braess edges form is asymptotically the same as that on $\mathcal{P}_{k}$ for sufficiently large $k$; and if $p$ dominates over $k$, then for sufficiently large $p$, edges added where $s<\frac{l}{l^2-l+1}$ (close to the single pendent vertex) tend to be Braess, while edges added where $s>\frac{l}{l^2-l+1}$ (close to the opposite side to the single vertex) do not. 

We now investigate asymptotic behaviours of $B_3$.

\begin{lemma}\label{lem:beta_0}
	Let $s_1(l,k,p)$ be the expression defined in \eqref{s_11}. There exists $L_k>0$ depending on $k$ between $\sqrt{k}$ and $\sqrt{3k}$ such that $s_1(l,k,p)>0$ for $3\leq l < L_k $, and $s_1(l,k,p)\le 0$ for $L_k \leq  l\leq k$. 
	Moreover, we have 
	$$\sum_{l=3}^{\lceil L_k\rceil-1} \lceil s_1(l,k,p)\rceil\leq B_3\leq 2\sum_{l=3}^{\lceil L_k\rceil-1} \lceil s_1(l,k,p)\rceil.$$
\end{lemma}
\begin{proof}
	For fixed \(k\) and $p$, the number of Braess edges on $\mathcal{P}_k$ for $\mathcal{B}_{k,p}$ is equal to the number of integers $l\ge 3$ satisfying \eqref{range of s}, that is,
	\begin{equation}\label{eq:count (45)}
		\max\{0,\lceil s_1(l,k,p) \rceil \}+ \max\{0,(k-l+1 - \lfloor s_2(l,k,p)\rfloor )\}.   
	\end{equation}
	Since $(k+p-l) -s_2(l,k,p) = s_1(l,k,p)$ by definition, we have $(k-l+1)-s_2(l,k,p)\leq s_1(l,k,p)$, so the inequalities
	\begin{equation}\label{eq:ineq for s1,s2}
		0 \le \max\{0,(k-l+1 - \lfloor s_2(l,k,p)\rfloor )\} \le \max\{0,\lceil s_1(l,k,p)\rceil\} 
	\end{equation}
	hold. Then using \eqref{eq:count (45)} and \eqref{eq:ineq for s1,s2}, we have
	$$\sum_{l\in L} \lceil s_1(l,k,p)\rceil\leq B_3\leq 2\sum_{l\in L} \lceil s_1(l,k,p)\rceil,$$
	where $L$ is the set of integers $l$ such that $s_1(l,k,p) >0$ for fixed $k$ and $p$. Note that $s_1(l,k,p)> 0$ if and only if $(k+p-l)^2-4C(l,k,p)> 0$. One can compute that
	\[
	(k+p-l)^2-4C(l,k,p) = ~\frac{N(l,k,p)}{3(k + p - 1)(l^2-l+1)},
	\]
	where
	\begin{multline*}
		N(l,k,p) = 3(k + p - 1)l^3 - 4(k+p-1)(k+p)l^2 +4(3k - 2)p^2 + (12k^2 - 16k + 8)p + 4k(k-1)^2.
	\end{multline*}
	So, $s_1(l,k,p)> 0$ if and only if $N(l,k,p)>0$. 
	We can find that for $k\geq 4$,
	\begin{align*}
		&N(\sqrt{k},k,p) = 3k^{5/2} + 4(p - 1)k^2 + 3(p-1)k^{\frac{3}{2}} + 4(2p-1)(p-1)k -8p(p-1) >0,\\
		&N(\sqrt{3k},k,p) = - 8k^3 + 9\sqrt{3}k^{\frac{5}{2}} - 4(3p-1)k^2 +9\sqrt{3}(p-1)k^{\frac{3}{2}} - 4(p-1)k - 8p(p-1) <0,\\
		&N(k,k,p) =  -k^4 - 5(p-1)k^3  -4(p^2 -4p +2)k^2 + 4(3p-1)(p-1)k - 8p(p-1)<0.
	\end{align*}
	We see that $N(l,k,p)$ is a cubic polynomial in $l$, and $N(0,k,p)>0$. Moreover,  the coefficients of $l^3,l^2$, and $l$ are positive, negative, and $0$, respectively. 
	It follows that given $k\geq 4$, there exists $L_k>0$ between $\sqrt{k} $ and $\sqrt{3k}$ such that $N(l,k,p)>0$ for $3\leq l < L_k$; and $N(l,k,p) \le 0$ for $L_k\leq  l\leq k$. (See Figure~\ref{fig:cubicPoly} for the illustration of $N(l,k,p)$ with respect to $l$.)
	This completes the proof.
	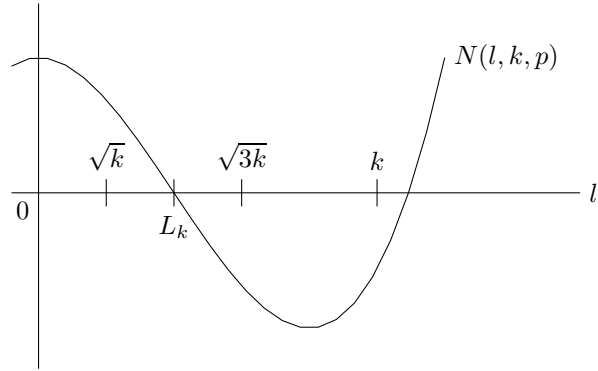
\begin{figure}[h!]
		\centering
		\begin{tikzpicture}[domain=-0.2:3, scale = 1.8]
			\draw plot (\x,{ (0.5)*(\x)^3 - (1.5)*(\x)^2+1}) node[right] {$N(l,k,p)$} ;
			\draw (-0.2,0) -- (4,0) node[right] {$l$};
			\draw (0,-1.3) -- (0,1.4);
			\draw (1,0.1) -- (1,-0.1) node[below]{ $L_k$}  ;
			\draw (0,0) node[below left]{ $0$} ;
			\draw (0.5,0.1) node[above]{ $\sqrt{k}$} -- (0.5,-0.1);
			\draw (1.5,0.1) node[above]{ $\sqrt{3k}$} -- (1.5,-0.1) ;
			\draw (2.5,0.1) node[above]{ $k$} -- (2.5,-0.1) ;
		\end{tikzpicture}\label{fig:cubicPoly}\caption{An illustration of $N(l,k,p)$.}
	\end{figure}
\end{proof}

From \eqref{eqn:c(l,k,p)}, we can recast $C(l,k,p)$ as follows:
\begin{align*}
	C(l,k,p) = \frac{(k+p)^2}{4}\left(H(l)+G(l)\right),
\end{align*}
where
\begin{align*}
	H(l) = &~1-\frac{kp}{(k+p)^2}\frac{4l}{l^2-l+1}-\frac{k}{k+p-1}\frac{k^2}{(k+p)^2}\frac{4l}{3(l^2-l+1)},\\
	G(l) =&~ \left(-2l+\frac{4l(l^2-1)}{3(l^2-l+1)}+\frac{4l}{l^2-l+1}\frac{p}{k+p}+\frac{k}{k+p-1}\frac{4l}{l^2-l+1}\frac{k}{k+p}\right)\frac{1}{k+p}\\
	&~+\left(-\frac{l(l^2-l)}{l^2-l+1}-\frac{k}{k+p-1}\frac{8l}{3(l^2-l+1)}\right)\frac{1}{(k+p)^2}.
\end{align*}
Then we see that
\[   \frac{s_1(l,k,p)}{k+p}  \;=\;   \frac{1}{2}\left( 1-\frac{l}{k+p}  -
\sqrt{H(l) + G(l)}      \right).
\]

\begin{theorem}\label{thm:B3}
	Suppose that $\frac{p}{k}\rightarrow \beta \in [0,\infty]$ as $k\rightarrow\infty$. Then \[B_3 = 
	\Theta(k \ln (k)).
	\]
\end{theorem}
\begin{proof}
	Let $s_1(l,k,p)$ be the expression defined in \eqref{s_11}, and let
	\begin{align*}
		s_1^*(k)   \;=\;   \sum_{l=3}^{\lceil L_k\rceil-1}  \lceil  s_1(l,k,p) \rceil \, ,
	\end{align*}   
	where $L_k$ is the positive number between $\sqrt{k}$ and $\sqrt{3k}$ defined in Lemma~\ref{lem:beta_0}. We write
	\begin{align}\label{s_1}
		\frac{s_1^*(k)}{k+p}  \;=\;   \frac{1}{2}\left(P_0(k)+P_1(k)+P_2(k)\right),
	\end{align}    
	where
	\begin{align*}
		P_0(k)  \; & = \;   \frac{1}{k+p}  \sum_{l=3}^{\lceil L_k\rceil-1} 
		\left( \lceil s_1(l,k,p) \rceil  -   s_1(l,k,p) \right),
		\\
		P_1(k)  \;& = \;    \sum_{l=3}^{\lceil L_k\rceil-1}  \left(  1-\sqrt{H(l)} \right),
		\\
		P_2(k) \; & = \;   \sum_{l=3}^{\lceil L_k\rceil-1}  \left(  -\frac{l}{k+p}  
		+\sqrt{H(l)}-\sqrt{H(l)-G(l) }\right)\,.
	\end{align*}
	
	We first show that $P_0(k)+P_2(k)$ is $o(\ln k)$.   As done in the proof of Proposition~\ref{prop-asymp}, we have $P_0(k)=O(\sqrt{k}/k)$.
	For $P_2(k)$, we use the identity (\ref{eq.diffsquares}) to obtain (for some constant $D$)
	\begin{align*} 
		\left| \sqrt{H(l)}-\sqrt{H(l)-G(l) }  \right|
		\; 
		=\;  \frac{|G(l)|}{   \sqrt{H(l)}+\sqrt{H(l)-G(l) } }\;\leq \;\frac{|G(l)|}{\sqrt{H(l)}} \;\leq \;  \frac{  1 }{\sqrt{H(l)}}\frac{Dl}{(k+p)}  \,.
	\end{align*}
	Since $\lim_{l\rightarrow \infty}H(l)=1$ and    $l\leq L_k \leq \sqrt{3k}$, 
	it follows that $P_2(k)$ is bounded by a constant.
	
	Now we shall examine $P_1(k)$. Let $$H_1(l) =1 - \frac{(k+p)^2}{(k+p)^2}\frac{4l}{l^2-l+1}-\frac{k+p-1}{k+p-1}\frac{(k+p)^2}{(k+p)^2}\frac{4l}{3(l^2-l+1)}.$$ 
	Since $1 - H(l)<1-H_1(l)$, we see that
	\begin{align*}
		P_1(k)\leq \sum_{l=3}^{\lceil L_k\rceil-1}  \left(  1-\sqrt{H_1(l)} \right).
	\end{align*}
	Using the identity
	\begin{equation}
		\label{eq.diffsquares}
		\sqrt{x}-\sqrt{y}  \;=\;  \frac{x-y}{\sqrt{x}+\sqrt{y}} \,,
	\end{equation}
	we observe that 
	\[
	\lim_{l\rightarrow\infty}\frac{1-\sqrt{H_1(l)}}{1/l}  \;=\;
	\lim_{l\rightarrow\infty}\frac{l(1-H_1(l))}{1+\sqrt{H_1(l)}} \; = \; \frac{8}{3}.      
	\]      
	Now we use an elementary exercise in analysis:  if $\{a_{l}\}$ and $\{b_{l}\}$ are two nonnegative real sequences such that $\lim_{l\rightarrow\infty}a_{l}/b_{l}=A>0$, and
	if $\sum_{l=1}^{\infty}b_{l}$ diverges,  then
	\begin{equation*}
		\lim_{n\rightarrow\infty}\frac{ \sum_{l=1}^{n}a_{l} }{ \sum_{l=1}^{n}b_{l}}  \;=\; A\,.
	\end{equation*}
	Since $\sum_{l=1}^{n}1/l\,\sim \,\ln n$ and $\sqrt{k} < L_k < \sqrt{3k}$, it follows that 
	\[    
	P_1(k)=O(\ln (k)).
	\]
	
	Suppose $k\geq \beta p$ for some $\beta>0$. Then 
	$$1 - H(l)\;> \;\frac{k}{k+p}\frac{k^2}{(k+p)^2}\frac{4l}{3(l^2-l+1)}\geq  \left(\frac{\beta}{1+\beta}\right)^3\frac{4l}{3(l^2-l+1)}.$$
	As done for obtaining $O(\ln(k))$, it follows that $$P_1(k)=\Omega(\ln k).$$
	
	Finally, consider the case for $p>k$. We see that $$1-\sqrt{H(l)}  = \frac{1-H(l)}{1+\sqrt{H(l)}} > \frac{1-H(l)}{2}\;>\; \frac{p}{k+p}\frac{k}{k+p}\frac{2l}{l^2-l+1}\;>\; \frac{1}{2}\frac{k}{k+p}\frac{2}{l}=\frac{k}{k+p}\frac{1}{l}.$$
	Hence, $$P_1(k) = \Omega(k\ln(k)/(k+p)).$$
	From \eqref{s_1}, $$s_1^*(k) = \Omega(k \ln (k)).$$
	
	The theorem follows.
\end{proof}

\end{document}